\def\VERSION{20.8.2017}  
\def\users{us}\def\friends{no} 
\def\users{world} 
\documentclass[11pt]{article} 
\RequirePackage{fix-cm}

\usepackage{color}
\usepackage{amsmath}
\usepackage{amsthm}
\usepackage{amssymb}
\usepackage{epsfig}
\usepackage{psfrag}
\usepackage{graphicx}
\usepackage{textcomp}
\usepackage{pgf}
\numberwithin{equation}{section}
\usepackage{upgreek} 
\newtheorem{theorem}{Theorem}[section]

\newtheorem{definition}[theorem]{Definition}

\newtheorem{proposition}[theorem]{Proposition}
\newtheorem{corollary}[theorem]{Corollary}
\newtheorem{remark}[theorem]{Remark}

\usepackage{mathrsfs,cite}
\usepackage{ifthen}
\usepackage{ulem}
\usepackage{cancel}\ifthenelse{\equal{\users}{world}}
{
\newcommand{\REM}[1]{}

	\newcommand{\DELETE}[1]{}

        \newcommand{\COMMENT}[1]{}
        \newcommand{\TCOMMENT}[1]{}

}	
{

\definecolor{brown}{rgb}{0.6,0.2,0.2}
\newcommand{\REM}[1]{\marginpar{\bfseries\tiny{\color{blue}#1}}}

 \newcommand{\COMMENT}[1]{{\color{blue}\uuline{#1}\color{black}}} 
 \newcommand{\DELETE}[1]{{\color{brown}\cancel{#1}\color{black}}}

 \newcommand{\TCOMMENT}[1]{{\color{blue}{ #1}}}

\newcount\hour \newcount\minute
\hour=\time
\divide \hour by 60
\minute=\time
\loop \ifnum \minute > 59 \advance \minute by -60 \repeat
}

\ifthenelse{\equal{\friends}{for-friends}}
{
}

\newcommand\DT[1]{\mathchoice
                 {{\buildrel{\hspace*{.1em}\text{\LARGE.}}\over{#1}}}
                 {{\buildrel{\hspace*{.1em}\text{\Large.}}\over{#1}}}
                 {{\buildrel{\hspace*{.1em}\text{\large.}}\over{#1}}}
                 {{\buildrel{\hspace*{.1em}\text{\large.}}\over{#1}}}}
\newcommand\DDT[1]{\mathchoice
   {{\buildrel{\hspace*{.13em}\text{\LARGE.\hspace*{-.13em}.}}\over{#1}}}
   {{\buildrel{\hspace*{.1em}\text{\Large.\hspace*{-.1em}.}}\over{#1}}}
   {{\buildrel{\hspace*{.1em}\text{\large.\hspace*{-.1em}.}}\over{#1}}}
   {{\buildrel{\hspace*{.1em}\text{\large.\hspace*{-.1em}.}}\over{#1}}}}
 \newcommand\DDDT[1]{\mathchoice
   {{\buildrel{\hspace*{.1em}\text{\Large.\hspace*{-.1em}\Large.\hspace*{-.1em}.}}\over{#1}}}
   {{\buildrel{\hspace*{.1em}\text{\large.\hspace*{-.1em}\large.\hspace*{-.1em}.}}\over{#1}}}
   {{\buildrel{\hspace*{.1em}\text{\large.\hspace*{-.1em}\large.\hspace*{-.1em}.}}\over{#1}}}
   {{\buildrel{\hspace*{.1em}\text{\large.\hspace*{-.1em}\large.\hspace*{-.1em}.}}\over{#1}}}}

\newcommand{\Vdots}{\mathchoice{\:\begin{minipage}[c]{.1em}\vspace*{-.4em}$^{\vdots}$\end{minipage}\;}{\:\begin{minipage}[c]{.1em}\vspace*{-.4em}$^{\vdots}$\end{minipage}\;}{\:\tiny\vdots\:}{\:\tiny\vdots\:}}
\renewcommand{\d}{{\rm d}}
\newcommand{\eq}{\eqref}

\newcommand{\eps}{\varepsilon}

\newcommand{\calG}{\mathcal G}

\newcommand{\R}{\mathbb{R}}
\newcommand{\N}{\mathbb{N}}
\newcommand{\bbI}{\mathbb{I}}

\newcommand{\bbC}{\mathbb{C}}
\newcommand{\bbD}{\mathbb{D}}

\newcommand{\pl}{\partial}

\newcommand{\OF}{\Omega_{_{\rm F}}}
\newcommand{\OS}{\Omega_{_{\rm S}}}
\newcommand{\QF}{Q_{_{\rm F}}}
\newcommand{\QS}{Q_{_{\rm S}}}
\newcommand{\MAT}{_{_{\rm M}}}

\newcommand{\Cdot}{\hspace{-.1em}\cdot\hspace{-.1em}}
\newcommand{\Colon}{\hspace{-.15em}:\hspace{-.15em}}
\newcommand{\In}{\!\in\!}
\newcommand{\dev}{\mathop{\mathrm{dev}}}
\newcommand{\Rsym}{\mathbb R^{3\times 3}_{\rm sym}}
\newcommand{\Rdev}{\mathbb R^{3\times 3}_{\rm dev}}
\newcommand{\Rsph}{\mathbb R^{3\times 3}_{\rm sph}}
\renewcommand\r{\varrho}
\newcommand{\Dir}{{\scriptscriptstyle\mathrm{D}}}  

\begin{document}
\begin{sloppypar}


\noindent{\huge\bf
Seismic waves and earthquakes in\\a global monolithic model.
}\bigskip\bigskip

\noindent{\Large Tom\'a\v s Roub\'\i\v cek$^{1,2}$}\bigskip\bigskip

%
%
\noindent{$\,\!^1$ Mathematical Institute, Charles University,
Sokolovsk\'a 83, CZ-18675~Praha~8, Czech Republic\\
$^2$ 
Institute of Thermomechanics, CAS,
Dolej\v skova 5, 
CZ-18200~Praha~8, Czech Republic 
}

\bigskip

{\small
\noindent{\it Abstract}:
The philosophy that a single ``monolithic'' model can ``asymptotically'' 
replace and couple in a simple elegant way several specialized models relevant 
on various Earth layers is presented and, in special situations, also 
rigorously justified. In particular, global seismicity and tectonics is 
coupled to capture e.g.\ (here by a simplified model) ruptures of 
lithospheric faults generating seismic waves which then propagate through the 
solid-like mantle and inner core both as shear (S) or pressure (P) waves, 
while S-waves are suppressed in the fluidic outer core and also in the oceans. 
The ``monolithic-type'' models have the capacity to describe all the 
mentioned features globally in a unified way together with corresponding 
interfacial conditions implicitly involved, only when scaling its parameters 
appropriately in different Earth's layers. Coupling of seismic waves with 
seismic sources due to tectonic events is thus an automatic side effect. 
The global ansatz is here based, rather for an illustration, only on a 
relatively simple Jeffreys' viscoelastic damageable material at small strains 
whose various scaling (limits) can lead to Boger's viscoelastic fluid or even 
to purely elastic (inviscid) fluid. Self-induced gravity field, Coriolis, 
centrifugal, and tidal forces are counted in our global model, as well. 
The rigorous mathematical analysis as far as the existence of
solutions, convergence of the mentioned scalings, and energy conservation 
is briefly presented.

\bigskip


\bigskip


} 

\bigskip

\section{Introduction}

Global geophysical models typically (have to) deal with several very 
different phenomena and couple various models due to the layered character 
of our planet Earth as well as our Moon and many other terrestrial planets 
or other moons. This paper wants to demonstrate (and, in special
situations, also rigorously justify) the philosophy that a single 
``monolithic'' model can replace several specialized models coupled together.
Such a single model can be  easier to implement on computers algorithmically. 
Of course, computationally, such a monolithic-type model may be not 
always easier to produce really relevant simulations at computers 
we have at our disposal nowadays. As purely seismic global 3D models 
are already treatable as well as their coupling with earthquake source 
at least locally, cf.\ e.g.\ 
\cite{KomTro02SESG,KomTro02SESG2,LayWal95MGS,ToCaMa09SSLV}
as well as \cite{HBAD09DERC,HuAmHel14ERMW,KaLaAm08SEMS,PPAB12TDDR}, 
respectively, there 
is a hope that the coupling monolithic approach may become more amenable 
in future within ever increasing computer efficiency.

The phenomena we have in mind in this paper involve {\it global seismicity} and 
{\it tectonics}. In particular the latter involves e.g.\ {\it ruptures of 
lithospheric faults} generating seismic waves which then propagate through the 
solid-like silicate mantle and iron-nickel inner core both in shear (S) or 
pressure (P) modes. In contrast to the P-waves (also called primary or
compressional waves), the S-waves (also called secondary waves) are suppressed 
in the fluidic iron-nickel outer core and also in the water oceans 
(where P-waves emitted from the earthquakes in the crust may manifest as 
Tsunami at the end). 

Typically, Maxwell-type rheologies are used in geophysical models
of solid 
mantle to capture long-term creep effects up to $10^5$ yrs. Sometimes, 
more attenuation of Kelvin-Voigt type is also involved, which leads 
to the {\it Jeffreys rheology} \cite{LyHaBZ11NLVE}, cf.\ also \cite{Roub17GMHF}.
This seems more realistic in particular because it covers (in the limit)
also the Kelvin-Voigt model applied to the volumetric 
strain whereas the pure Maxwell rheology allowing for big creep
during long geological periods is not a relevant effect in the volumetric part.
We therefore take the Jeffreys model as a basic global ``monolithic'' 
ansatz and, in various limits in the deviatoric and the volumetric parts, 
we model different parts of the planet Earth.

Respecting the solid parts of the model, we use the 
{\it Lagrangian description}, i.e.\ here all equations are formulated 
terms of displacements rather than velocities, while the reference
and the actual space configurations automatically coincide with each other in our
small strain (and small displacement) ansatz, which is well relevant
in geophysical short-time scales of seismic events. 

In the solid-like part, various inelastic processes are considered 
to model tectonic earthquakes on lithospheric faults together with
long healing periods in between them, as well as aseismic slips, and
various other phenomena. To this goal, many internal variables
may be involved as aging/damage, inelastic strain, porosity, water content, 
breakage, and temperature, cf.\ \cite{LyHaBZ11NLVE,LyaBZ14DBRM}.
On the other hand, those sophisticated models are focused on 
rather local events around the tectonic faults without ambitions 
to be directly coupled with the global seismicity. 
Cf.\ also \cite{Roub14NRSF,RoSoVo13MRLF} for models of 
rupturing lithospheric faults and, in particular, a relation with 
the popular Dieterich-Ruina rate-and-state friction model. Here, rather for 
the lucidity of the exposition, we reduce the set of internal variables 
to only one scalar variable, namely damage/aging, which however has 
a capacity to trigger a spontaneous rupture (so-called {\it dynamic triggering}) 
with emission of seismic waves and,
in a certain simplification, can serve as a seismic source coupled 
with the overall global model. Also, this simple model already 
will well illustrate mathematical difficulties related to nonlinearities 
in the solid parts coupled with linear but possibly hyperbolic fluidic regions.

Let us emphasize that usual models are focused only either on propagation of 
seismic waves along the whole globe while their source is considered given, or 
on description of seismic sources due to tectonic events, 
but not their mutual coupling.
If a coupling is considered, then is concerns rather local models 
not considering the layered structure of the whole planet, cf.\ e.g.\  
\cite{BeZi01DRRM,BZamp09SRRS,HuAmHel14ERMW,KaLaAm08SEMS,LHAB09NDRW}.
The reality ultimately 
captures very different mechanical properties of different 
layers of the Earth, in particular the mantle and the inner core which are 
solid from the short-time scales versus the outer core and the oceans 
which are fluidic even on the short-time scales.
 
The goals of this article are:
\begin{itemize}
\vspace*{-.2em}\item[$\upalpha$)\!]to propose a model that might
capture simultaneously propagation of seismic waves over the whole planet
and their nonlinearly behaved sources (like ruptures of tectonic faults, 
here modelled only in a very simplified way for a relatively lucid 
illustration of the model procedure), both mutually coupled. 
\vspace*{-.0em}\item[$\upbeta$)\!]
by proper scaling to approximate viscoelastic Boger-type \cite{Bog77HECV}
fluids that is relevant in outer core and in the oceans (with a very low 
viscosity) where S-waves can then only slightly penetrate the outer core or
the oceans but are fast attenuated, while P-waves are only refracted.
\vspace*{-.0em}\item[$\upgamma$)\!]
by limiting further the viscosity in the outer core or the oceans to zero, 
further to approximate this viscoelastic fluid towards  elastic (completely
inviscid) fluid, 
respecting the phenomenon that 
S-waves cannot penetrate into these fluidic regions and are fully reflected 
on the interfaces between the outer core and the mantle 
(=Gutenberg's discontinuity) and the inner core
or on the ocean beds, i.e.\ between $\OS$ and $\OF$, while P-waves 
propagate through these interfaces, being both refracted and reflected on them. 
\vspace*{-.0em}\item[$\updelta$)\!]
perform the rigorous analysis as far as the existence of solutions,
a-priori estimates in specific norms, and convergence towards other models
that justifies the particular models, their energetics and asymptotics,
and can support numerical stability and convergence when discretised
and implemented on computers.  
\end{itemize}
\begin{figure}
\begin{center}
\psfrag{OL}{\small $\OF$}
\psfrag{OS}{\small $\OS$}
\psfrag{1200}{\small $1200\,$km}
\psfrag{3480}{\small $3480\,$km}
\psfrag{6371}{\small $6371\,$km}
\psfrag{ocean}{\scriptsize\bf\begin{minipage}[t]{15em}\hspace*{-2em}OCEANS
 - elastic fluid ($v_{_{\rm P}}\sim1.5\,$km/s, $v_{_{\rm S}}=0$,
very low viscosity $10^{-3}\,$Pa\,s,
bulk modulus $\sim$2GPa.)\end{minipage}}
\psfrag{mantle}
{\scriptsize\bf\begin{minipage}[t]{12em}\hspace*{-2em}CRUST \& MANTLE - solid 
(density
$\sim$ 5000 kg/m$^3$,
viscosity $10^{21}-10^{23}\,$Pa\,s, $v_{_{\rm P}}=8-13\,$km/s,
$v_{_{\rm S}}=5-7\,$km/s)
\end{minipage}}
\psfrag{outer core}{\scriptsize\bf\begin{minipage}[t]{13em}\hspace*{-2em}OUTER CORE - elastic fluid 
($v_{_{\rm P}}=8-10\,$km/s, $v_s=0$,
density $\sim$ 11000 kg/m$^3$,
very low viscosity $10^{-2}-10^{4}\,$Pa\,s)\end{minipage}}
\psfrag{inner core}{\scriptsize\bf\begin{minipage}[t]{14em}\hspace*{-2em}INNER CORE - solid
(high density $\sim$ 13000 kg/m$^3$, $v_{_{\rm P}}\sim11\,$km/s, $v_{_{\rm S}}\sim4\,$km/s,
viscosity $10^{14}-10^{15}\,$Pa\,s)\end{minipage}}
\hspace*{2em}\includegraphics[width=35em]{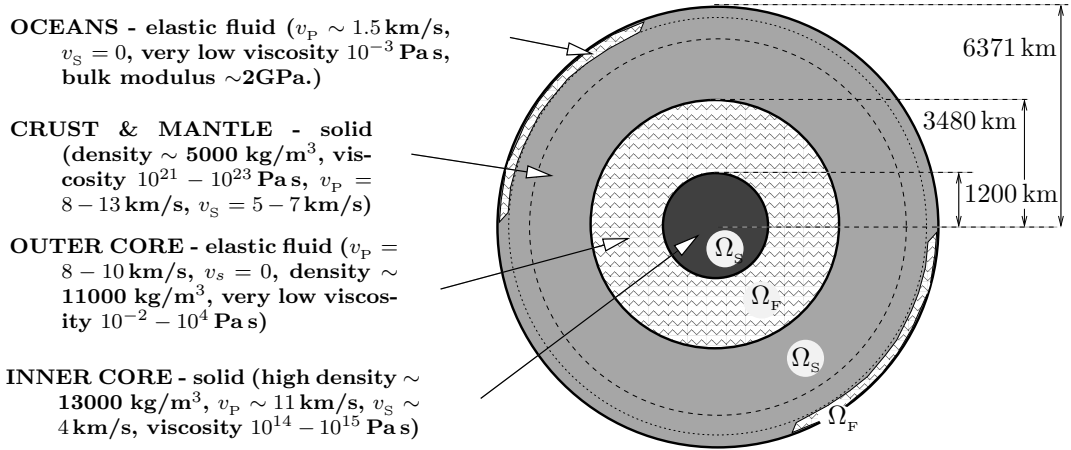}
\end{center}
\vspace*{-.6em}
\caption{\sl The very basic layered structure of our planet Earth with some 
(only very rounded) geometrical and rheological data and the notation for
the solid-type domains $\OS$ and the fluid-type domains $\OF$. Here, 
$v_{_{\rm P}}$ and 
$v_{_{\rm S}}$ stand for the velocity of the P- and S-waves, respectively, 
$v_{_{\rm S}}=0$ indicating that only the P- but not S-waves can propagate through
the particular layer.}
\label{fig-earth}
\end{figure}
The plan of this article is the following: In Section~\ref{sec2}, the general 
monolithic model is introduced and also specialized for viscoelastic 
linear isotropic material undergoing possibly damage in the elastic shear part.
In Sections~\ref{sec-Boger} and \ref{sec4}, we formulate the limit toward 
viscoelastic and purely elastic fluids, respectively. This limit concerns
only the sudbomain $\OF$, i.e.\ the outer core and oceans. To make the 
presentation more lucid and readable also for the geophysical specialist 
outside mathematical-analyst community, the mathematical analysis of the 
models and necessarily a bit technical proofs of the statements formulated in 
Sections~\ref{sec2}--\ref{sec4} are intentionally presented as Appendix later in 
Section~\ref{sec-anal}. For it, we use a conceptually constructive approximation
by the Galerkin method, which is in some variant also used in geophysical 
literature, cf. e.g.\ \cite{PPAB12TDDR}.

\section{Philosophy and an example of a monolithic model}\label{sec2}

\def\ein{\pi}

Rather for illustration of our main focus to coupling of solid and fluidic
models, we ignore most of the above mentioned internal variables and 
keep only the scalar-value aging/damage variable, denoted by $\alpha$,
and then we consider the mentioned Jeffreys rheology combined with damage 
influencing only the elastic but not viscous part. 
The differential relation governing this rheology is of the type 
``$\,\DT\sigma{+}\sigma=\DDT e{+}\DT e\,$'' with some specific 
coefficients (here ignored), some of them being possibly subject to damage. 
Here and in what follows, the dot-notation $(\cdot)\!\DT{^{}}$
stands for the time derivative $\frac{\partial}{\partial t}$.
Alternatively, to make the damage dependence more 
lucid, one can use system of 
first-order equations in time when implementing the concept of
internal variables, i.e.\ here introducing a Maxwellian-type 
creep strain $\ein$ 
and making a Green-Naghdi \cite{GreNag65GTEP}
additive decomposition of the total strain 
\begin{align}\label{e-decompose}
e(u)=e_{\rm el}+\ein.
\end{align}
So, considering also the damage/aging variable, altogether our set of 
internal variables will be $(\ein,\alpha)$.

The main ingredients of the model are then, beside the mass density 
$\varrho=\varrho(x)$, the specific stored energy and the dissipation 
potential as a property of the material 
\begin{subequations}\label{phi-zeta}
\begin{align}\label{phi}
&\varphi\MAT=\varphi\MAT(x,u,e_{\rm el},\alpha,\nabla\alpha,\phi,\nabla\phi),
\\
&\zeta\MAT=\zeta\MAT(x,\DT e_{\rm el},\DT\ein,\DT\alpha).
\end{align}\end{subequations}
Note that $\varphi\MAT$ and $\zeta\MAT$ naturally do not explicitly 
depend on the total strain $e(u)$ but 
rather on the elastic strain $e_{\rm el}$ and possibly on the 
creep strain $\ein$ and their rates; in fact, $\varphi\MAT$ is independent
of $\ein$ because naturally no hardening-like effects are relevant in
geophysical models. In what follows, we will often omit 
the explicit dependence on $x$ variable.

For readers' convenience, Table~1 summarizes the main nomenclature used 
in this paper.

\begin{table}[ht]
\centering
\fbox{
\begin{minipage}[t]{0.45\linewidth}
\normalsize

$u$ displacement, valued in $\R^3$,

$\ein$ creep strain, valued in $\Rsym$,

$\alpha$ damage/aging, valued in $[0,1]$,

$\phi$ potential of the gravity field,

$e(u)=\frac12\nabla u^\top\!{+}\frac12\nabla u$ total 
 strain,

$e_{\rm el}=e(u){-}\ein$ elastic strain,

$\varphi\MAT,\,\varphi$ stored energies,

$\zeta\MAT,\,\zeta$ potentials of dissipative forces,

$\eta$ potential of damage dissipation,

$f_{_{\rm COR}}(\DT u)=2\varrho\omega{\times}\DT u$ Coriolis force,

\smallskip

$\Rsym=\{A\In\R^{3\times3};\ A^\top=A\}$,

\smallskip

$\Rdev=\{A\In\Rsym;\ {\rm tr}\,A=0\}$,

\smallskip

$\Rsph=\{A\In\Rsym;\ A=a\bbI,\ a\In\R\}$,

\end{minipage}
\hfill 
\begin{minipage}[t]{0.52\linewidth}
\normalsize

$\varrho$ given mass density profile,

$\varrho_{\rm ext}^{}
$ external mass causing tidal forces,

$K_{_{\rm E}}$, $G_{_{\rm E}}$ bulk and shear elastic moduli,

$K_{_{\rm MX}}$, $G_{_{\rm MX}}$ 
Maxwell viscous moduli,

$K_{_{\rm KV}}$, $G_{_{\rm KV}}$ 
Kelvin-Voigt viscous moduli, 

$g$ the gravitational constant




$\omega$ given angular velocity of Earth rotation,

$\Omega\subset\R^3$ the reference domain (the Earth),

$\OF\subset\Omega$ the outer core 
and the oceans, 

$\OS\subset\Omega$ the mantle and the inner core,

$\Gamma$ the boundary of $\Omega$ (=Earth surface),

$I:=[0,T]$ the fixed time interval, 

$Q:=I\times\Omega$, $\ \QS:=I\times\OS$, $\ \QF:=I\times\OF$,

$\Sigma:=I\times\Gamma$.

\end{minipage}
}
\centerline{{\normalsize\bf Table~1.}$^{^{^{^{^{}}}}}$ 
 {\normalsize\sl Summary of the basic notation used through the paper.}}
\label{Tab_Notation}
\end{table}

\begin{figure}
\psfrag{e}{\small $e=e(u)$}
\psfrag{a}{\small $\alpha$}
\psfrag{s}{\small $\sigma$}
\psfrag{s1}{\small $\sigma_{\rm sph}$}
\psfrag{s2}{\small $\sigma_{\rm dev}$}
\psfrag{e1}{\small ${\rm dev}\,e$}
\psfrag{ee1}{\small ${\rm dev}\,e_{\rm el}$}
\psfrag{e2}{\small ${\rm sph}\,e$}
\psfrag{ee2}{\small ${\rm sph}\,e_{\rm el}$}
\psfrag{G}{\small $G_{_{\rm E}}\!=G_{_{\rm E}}\!(\alpha)$}
\psfrag{Gmx}{\small $G_{_{\rm MX}}$}
\psfrag{Gkv}{\small $G_{_{\rm KV}}$}
\psfrag{K}{\small $K_{_{\rm E}}$}
\psfrag{Kmx}{\small $K_{_{\rm MX}}$}
\psfrag{Kkv}{\small $K_{_{\rm KV}}$}
\begin{center}
\includegraphics[width=31em]{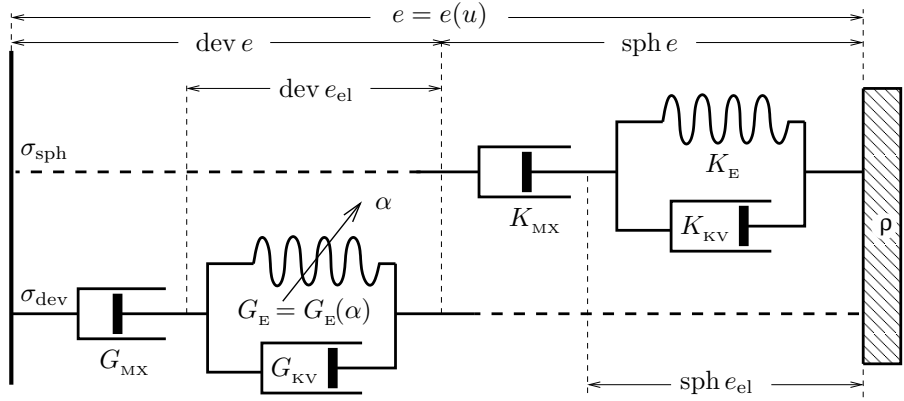}
\end{center}
\vspace*{-.6em}
\caption{\sl Schematic diagram for the viscoelastic Jeffreys rheology
which is subjected to damage $\alpha$ in the deviatoric part while
the spherical (volumetric) part is not subjected to damage.}
\label{fig-mixed-response}
\end{figure}

To formulate the equations for the dynamics of displacement $u$, it is more 
natural in the Kelvin-Voigt model (which is a part
of the Jeffreys rheology) to express the stored-energy and dissipation 
potentials 
rather in terms of the total strain $e=e(u)$ and the creep $\ein$ from 
\eq{e-decompose} as
\begin{subequations}\label{def-of-phi-zeta-}
\begin{align}
&\varphi(x,u,e,\ein,\alpha,\nabla\alpha,\phi,\nabla\phi)
=\varphi\MAT(x,u,e{-}\ein,\alpha,\nabla\alpha,\phi,\nabla\phi),
\\&\label{def-of-zeta-}
\zeta(\DT e,\DT \ein,\DT\alpha
)
=
\zeta\MAT(\DT e{-}\DT \ein,\DT \ein,\DT\alpha
).
\end{align}\end{subequations}
In terms of these potentials, our system in the abstract form 
is: 
\begin{subequations}\label{system}\begin{align}\label{system-u}
\varrho\DDT u+\varphi_u'(\nabla\phi)
-{\rm div}\,\big(
\zeta_{\DT e}'(e(\DT u),\DT \ein)+\varphi_e'(e(u),\ein,\alpha)\big)
&=f_{_{\rm COR}}(\DT u)
&&\text{in }\ \Omega,&&
\\
\zeta_{\DT \ein}'(e(\DT u),\DT \ein)+\varphi_\ein'(e(u),\ein,\alpha)
&=0
&&\text{in }\ \Omega,
\\
\zeta_\alpha'(\DT\alpha)
+\varphi_\alpha'(e(u),\ein,\alpha)
-{\rm div}(
\varphi_{\nabla\alpha}'(\nabla\alpha))&=0&&\text{in }\ \Omega,
\\
\varphi_\phi'(u,\phi)-{\rm div}(\varphi_{\nabla\phi}'(\nabla\phi))
&=\varrho_{\rm ext}(t)
&&\text{in }\ \R^3,
\end{align}\end{subequations}
where $(\cdot)'$ denotes the (partial) derivatives for which 
we already have assumed a certain particular form of \eq{phi-zeta}
so that the partial derivatives of $\varphi$ and $\zeta$ do not
depend on the full list of arguments, thus e.g.\ instead of 
$\varphi_e'(u,e(u),\ein,\alpha,\nabla\alpha,\phi,\nabla\phi)$ we wrote only 
$\varphi_e'(e(u),\ein,\alpha)$, etc. Note that, in \eq{system-u}, we 
used also a bulk force $f_{_{\rm COR}}(\DT u)$ which may not come from the
dissipation potential by a classical way
$f_{_{\rm COR}}(\DT u)=\zeta_{\DT u}$, which in particular
the case of the Coriolis (pseudo) force which itself does not 
vanish but nevertheless the dissipation rate $f_{_{\rm COR}}(\DT u)\Cdot\DT u=0$.

Rheological response under the volume variation
and under the shear may (and do) substantially differ from each other.
To distinguish these geometrical aspects,  the total strain is decomposed to
the {\it spherical} (also called hydrostatic or volumetric) 
{\it strain}\index{strain!spherical} 
and the {\it deviatoric} (also called shear) 
{\it strain}:\index{strain!deviatoric} 
\begin{align}\label{decomposition-of-strain}
e={\rm sph}\,e+{\rm dev}\,e\ \ \ \text{ with }\ \ \ 
{\rm sph}\,e=\frac{{\rm tr}\,e}3\bbI\ \ \ \text{ and }\ \ \ 
{\rm dev}\,e=e-\frac{{\rm tr}\,e}3\bbI,
\end{align}
where $\bbI$ is the identity matrix and ``tr'' denotes the trace of
a matrix. Note that the deviatoric and the spherical strains 
from \eq{decomposition-of-strain} are orthogonal to each other.
In terms of this decomposition, for example the isotropic elastic 
(Lam\'e) material at small strains has the quadratic stored energy
\begin{align}\nonumber
&\varphi_{_{\mbox{\scriptsize\rm Lam\'e}}}(e)=\frac12\lambda({\rm tr}\,e)^2+G_{_{\rm E}}|e|^2=
\frac12\Big(\lambda{+}\frac23G_{_{\rm E}}\Big)({\rm tr}\,e)^2
+G_{_{\rm E}}\Big|e-\frac{{\rm tr}\,e}3\bbI\Big|^2
\\&\nonumber\qquad\qquad=
\frac12\Big(3\lambda{+}2G_{_{\rm E}}\Big)|{\rm sph}\,e|^2
+G_{_{\rm E}}\big|{\rm dev}\,e\big|^2
\\&\qquad\qquad=\frac32K_{_{\rm E}}|{\rm sph}\,e|^2+G_{_{\rm E}}|{\rm dev}\,e|^2\ \ \ \text{ with }\ K_{_{\rm E}}=\lambda+\frac23G_{_{\rm E}}.
\label{Lame-K-G}\end{align}
where we used $|{\rm sph}\,e|^2=(\frac{{\rm tr}\,e}3)^2\bbI{:}\bbI=
(\frac{{\rm tr}\,e}3)^23=({\rm tr}\,e)^2/3$ and 
the mentioned orthogonality of ${\rm sph}\,e$ and ${\rm dev}\,e$,
where 
$K_{_{\rm E}}$ is called the bulk modulus and $\lambda$ is 
the (first) Lam\'e coefficient while $G_{_{\rm E}}$ is the shear modulus
(=\,the second Lam\'e coefficient).
Then the corresponding stress writes as:
\begin{align*}
\varphi_{_{\mbox{\scriptsize\rm Lam\'e}}}'(e)&=\lambda({\rm tr}\,e)
\bbI+2G_{_{\rm E}} e=:\sigma_{\rm dev}+\sigma_{\rm sph}
\\&\ \ \text{ with }\ \  
\sigma_{\rm dev}=2G_{_{\rm E}} {\rm dev}\,e\ \ \ \text{ and }\ \  
\sigma_{\rm sph}=3K_{_{\rm E}}{\rm sph}\,e.
\end{align*}

We will use this decomposition both for elastic and for the viscous parts
of our model and let the elastic bulk modulus $G_{_{\rm E}}$ and thus also 
$\varphi$ dependent on damage. Yet, this dependence would bring mathematical 
difficulties (like lack of the integrability of the term $\DT\alpha G_{_{E}}(\alpha)
|{\rm dev}\,e_{\rm el}|^2$ that would arise in \eq{system-alpha-w-}, for 
example). Here, various modifications of the Lam\'e-type model can help. 
One option is to consider higher strain gradients, i.e.\ the concept
of nonsimple continua, cf.\ also \cite{Roub17GMHF}. Other (conceptually
and technically simpler) option is to modify the stored energy by considering
an non-quadratic $e\mapsto\calG_{_{\rm E}}(\alpha,{\rm dev}\,e)$ instead
of $e\mapsto G_{_{E}}(\alpha)|{\rm dev}\,e|^2$. A canonical option one may have
in mind is
\begin{align}\label{KV-limit+}
\calG_{_{\rm E}}
(\alpha,e_{\rm el})=\calG_{_{\rm E}}
(\alpha,{\rm dev}\,e_{\rm el})=
\frac{G_{_{\rm E}}\!(\alpha)|{\rm dev}\,e_{\rm el}|^2}{\sqrt{1+|{\rm dev}\,e_{\rm el}|^2/E_{_{\rm M}}^2}}
\end{align}
for (presumably large) regularizing parameter $E_{_{\rm M}}$. 
In accord with Figure~\ref{fig-mixed-response}, in \eq{KV-limit+} and
in what follows, we quite naturally assume that $\calG_{_{\rm E}}$ does not depend 
on the spherical part of the strain which itself usually cannot make any
damage in the rock. The meaning of $E_{_{\rm M}}$ in \eqref{KV-limit+}
is clear from realizing that, if the deviatoric 
strain is substantially smaller in the sense 
$|{\rm dev}\,e_{\rm el}|<\!<E_{_{\rm M}}$, then the difference from the original 
Lame\'e-type model is negligible, i.e.\ the corresponding stress contribution 
$\partial_e 
\calG_{_{\rm E}}(\alpha,\,e_{\rm el})\sim 2G_{_{\rm E}}\!(\alpha)
{\rm dev}\,e_{\rm el}$.
On the other hand, this modification is effective for large strains and makes 
$e_{\rm el}\mapsto\varphi\MAT(x,u,e_{\rm el},\alpha,\nabla\alpha,\phi,\nabla\phi)$ 
nonquadratic but still keeps convex now with at most linear growth,
so that the elastic part of deviatoric stress is always below 
$G_{_{\rm E}}\!(\alpha)E_{_{\rm M}}$. Thus, for $E_{_{\rm M}}$ large, this is effectively
not a substantial (and physically well acceptable) modification 
of the original model. 

Thus we consider the {\it material stored energy} and the 
{\it material dissipation-force potential} from \eq{phi-zeta} governing the 
problem as 
\begin{subequations}\label{def-of-phi-zeta}
\begin{align}
&\varphi\MAT\begin{cases}=\varphi\MAT(x,u,e_{\rm el},\alpha,\nabla\alpha,\phi,\nabla\phi)=
\frac32K_{_{\rm E}}|{\rm sph}\,e_{\rm el}|^2+
\calG_{_{\rm E}}\!(\alpha,{\rm dev}\,e_{\rm el})
\hspace{-6em}
&
\\\qquad\qquad\qquad\qquad\qquad\quad+\frac{\kappa}2|\nabla\alpha|^2+\varrho\phi
+\varrho\nabla\phi{\cdot}u+\frac1{8\uppi g}|\nabla\phi|^2\hspace{-6em}&
\\\qquad\qquad\qquad\qquad\qquad\quad+\frac12\varrho\big(|\omega{\cdot}(x{+}u)|^2{-}|\omega|^2|x{+}u|^2\big)
&\text{if }x\in\Omega,
\\=
\varphi\MAT(\nabla\phi)=\frac1{8\uppi g}|\nabla\phi|^2&\text{if }x\in\R^3\!\setminus\!\Omega,
\end{cases}
\label{def-of-phi}\\
\nonumber
&\zeta\MAT=\zeta\MAT(x,\DT e_{\rm el},\DT \ein,\DT\alpha
)=\frac32K_{_{\rm KV}}|{\rm sph}\,\DT e_{\rm el}|^2
+G_{_{\rm KV}}|{\rm dev}\,\DT e_{\rm el}|^2
\\&\qquad\qquad\qquad\qquad\qquad\qquad
+\frac32K_{_{\rm MX}}|{\rm sph}\,\DT \ein|^2
+G_{_{\rm MX}}|{\rm dev}\,\DT \ein|^2
+\eta(\DT\alpha)
\label{def-of-zeta}
\end{align}\end{subequations}
with $g\doteq6.674\times10^{-11}$m$^3$kg$^{-1}$s$^{-2}$ the gravitational constant. 
For the $\phi$-terms in \eq{def-of-phi} see e.g.\ 
\cite[Formulae (16),(20)]{WooDeu09TOEF}. 
All the coefficients $K$'s, $G$'s, and $\kappa$'s 
are naturally considered defined on $\Omega$ and $x$-dependent (which is not 
explicitly written in \eq{def-of-phi-zeta} for brevity).
Also the mass density is considered  $x$-dependent
and, because of \eq{system4}, defined on the whole Universe $\R^3$ 
by putting $\varrho(x)=0$ if $x\in\R^3{\setminus}\Omega$.
Keeping $K_{_{\rm E}}$ constant in \eq{def-of-phi} reflects the phenomenon that
compression does not cause damage and the pure tension is not 
a relevant mode in geophysical models, so only shear does cause
damage because of dependence of the elastic shear modulus 
$G_{_{\rm E}}$ in \eq{KV-limit+} and \eq{def-of-phi} on $\alpha$.
The last term in \eq{def-of-phi} is the potential of the centrifugal
force $\varrho\omega{\times}(\omega{\times}(x{+}u))$
and, noteworthy, it violates coercivity due to the
term ${-}\frac12\varrho|\omega|^2|x{+}u|^2$, which reflects the real phenomenon
that, given the angular velocity $\omega$, the centrifugal force
can indeed inflate the body in an unlimited way in the direction orthogonal
to the rotation axis. Of course, in reality, this is either not 
observed in planetary/satellite bodies or the angular velocity cannot be taken
apriori fixed in some large-type bodies. 

For shorter notation, we define the 4th-order tensors corresponding
to the isotropic material expressed in terms of the $(K,G)$-moduli in 
\eq{def-of-phi-zeta}, namely
\begin{subequations}\label{C-D}
\begin{align}
&[\bbC(\alpha,e)]_{ijkl}=\,K_{_{\rm E}}\delta_{ij}\delta_{kl}
+\partial_{ee}^{}\calG_{_{\rm E}}\!(\alpha,e)(\delta_{ik}\delta_{jl}{+}\delta_{il}\delta_{jk}{-}2\delta_{ij}\delta_{kl}/3),
\\&[\bbD_{_{\rm MX}}]_{ijkl}=K_{_{\rm MX}}\delta_{ij}\delta_{kl}
+G_{_{\rm MX}}(\delta_{ik}\delta_{jl}+\delta_{il}\delta_{jk}-2\delta_{ij}\delta_{kl}/3),
\\&[\bbD_{_{\rm KV}}]_{ijkl}=K_{_{\rm KV}}\delta_{ij}\delta_{kl}
+G_{_{\rm KV}}(\delta_{ik}\delta_{jl}+\delta_{il}\delta_{jk}-2\delta_{ij}\delta_{kl}/3),
\end{align}\end{subequations}
where $\delta$'s denote the Kronecker symbol.
In what follows, we will take specifically 
$f_{_{\rm COR}}(\DT u)=2\varrho\omega{\times}\DT u$
which is the standard {\it Coriolis force} with
$\omega$ the given vector (angular velocity
related with the constant rotation of the Earth with respect
to the inertial system).
The general system \eq{system} then takes the more specific form:  
\begin{subequations}\label{system+}\begin{align}\nonumber
&\varrho\DDT u-{\rm div}\,\sigma
+f=0\ \ \text{ with }\ \sigma=\bbD_{_{\rm KV}}\DT e_{\rm el}
+\bbC(\alpha,{\rm dev}\,e_{\rm el})e_{\rm el}\ 
\\&\label{system-u+}
\hspace*{9.3em}\text{and }\ 
f=\varrho\big(\omega{\times}(\omega{\times}(x{+}u))+2\omega{\times}\DT u+\nabla\phi\big)&&\text{in }\ \Omega,
\\&\label{system-pi}
\bbD_{_{\rm MX}}\DT \ein
=\sigma
&&\text{in }\ \Omega,
\\&\label{system-alpha}
\pl\eta(\DT\alpha)
+
\partial_\alpha\calG_{_{\rm E}}(\alpha,{\rm dev}\,e_{\rm el})
-{\rm div}(\kappa\nabla\alpha)\ni0
&&\text{in }\ \Omega,
\\&\label{system4}
{\rm div}\Big(\frac1{4\uppi g}\nabla\phi+\varrho u\Big)=\varrho
+\varrho_{\rm ext}(t)
&&\text{in }\,\R^3.\!
\end{align}\end{subequations}
The bulk force $f$ in \eq{system-u} involves the centrifugal force, the 
Coriolis force, and the 
gravity force due to the self-induced gravity field. This last force plays 
a certain role in ultra-low frequency (i.e.\ very long wavelength) seismic 
waves, cf.\ e.g.\ \cite{ToCaMa09SSLV}. 
For such right-hand side in \eq{system-u+} together with 
\eq{system4} see e.g.\ \cite[Formulae (19),(22)]{WooDeu09TOEF}. Actually,
this self-gravity interaction results after a certain linearization 
of the system originally written at large strains, cf.\ also 
\cite{Braz17EGEG,BrHoHo17VFEE}. The external time-varying 
mass density $\varrho_{\rm ext}^{}=\varrho_{\rm ext}^{}(x,t)$ occurring in 
\eq{system4} allows for involvement of {\it tidal forces} arising from 
other bodies (stars, planets, moons) moving around.
The notation $\pl\eta$ in \eq{system-alpha} stands for a subdifferential
of the convex function $\eta$ which standardly generalizes the usual
derivative if $\eta$ is non-smooth and, by definition, the inclusion 
``$\ni$'' there actually means, for any $v\In\R$ and a.a.\ $(t,x)\in Q$, the inequality 
\begin{align}\label{system-alpha-w-}
(v{-}\DT\alpha(t,x))
\Big(\partial_\alpha\calG_{_{\rm E}}(x,\alpha(t,x),{\rm dev}\,e_{\rm el}(t,x))
-{\rm div}(\kappa\nabla\alpha(t,x))\Big)
+\eta(v)\ge\eta(\DT\alpha(t,x)).
\end{align}

To facilitate spontaneous rupture on some weak narrow layers (pre-existing 
faults in the solid crust) that leads to earthquakes, some 
``enough pronounced'' nonconvexity of $\varphi$ is necessary. One option was 
devised in \cite{LyaMya84BECS}, introducing non-convexity in terms
of $e_{\rm el}$ if damage develops enough, and used in a series of papers, 
see e.g.\ \cite{LyHaBZ11NLVE,LyaBZ14DBRM} and references therein. 
Mathematically rigorous setting of this model requires however 
some higher-order strain gradients (the so-called 2nd-grade nonsimple
materials) and is rather technical, cf.\ \cite{Roub17GMHF}. 
Alternative nonconvexity can be implemented through 
$\calG_{_{\rm E}}\!(\cdot,e)$ causing 
possible weakening when damage develops, being opposite to so-called cohesive 
damage. To avoid the mentioned mathematical technicalities, we have
chosen the latter option.

To see the energetics behind the system \eq{system}, we test the particular 
equations in \eq{system} by $\DT u$, $\DT\pi$, $\DT\alpha$, and $\DT\phi$, 
respectively. 
Thus we obtain (at least formally) the energy balance 
\begin{align}\nonumber
\frac{\d}{\d t}\bigg(\int_\Omega\frac\varrho2|\DT u|^2\,\d x+
\int_{\R^3}\varphi(u,e(u),\ein,\alpha,\nabla\alpha,\phi,\nabla\phi)\,\d x\bigg)\qquad\quad
\\[-.3em]+\int_\Omega\!\xi(e(\DT u),\DT\ein,\DT\alpha)\,\d x
=\int_{\R^3}\!\!\varrho_{\rm ext}(t)\DT\phi\,\d x\,,
\label{engr-abstract}\end{align}
where $\xi=\xi(\DT e,\DT\ein,\DT\alpha)$ is the dissipation rate related 
with the dissipation potential $\zeta=\zeta(\DT e,\DT\ein,\DT\alpha)$ by 
$\xi(\DT e,\DT\ein,\DT\alpha)
=(\DT e,\DT\ein,\DT\alpha)^\top\partial\xi(\DT e,\DT\ein,\DT\alpha)$.
Using \eq{def-of-phi-zeta} and \eq{C-D}, we can make it 
more specific for the system \eq{system+} when using also the calculus
\begin{align}\label{calculus}
\int_\Omega\varrho\nabla\phi\Cdot\DT u-{\rm div}(\varrho u)\DT\phi\,\d x
=\int_\Omega\varrho\nabla\phi\Cdot\DT u+\varrho u\Cdot\nabla\DT\phi\,\d x
=\frac{\d}{\d t}\int_\Omega\varrho\nabla\phi\Cdot u\,\d x.
\end{align} 
and realizing that, in view of  \eq{def-of-zeta-} and \eq{def-of-zeta},
\begin{align}\nonumber
\xi(\DT e,\DT\ein,\DT\alpha)&=3K_{_{\rm KV}}|{\rm sph}(\DT e{-}\DT\ein)|^2
+2G_{_{\rm KV}}|{\rm dev}(\DT e{-}\DT\ein)|^2
\\&\nonumber
\quad+3K_{_{\rm MX}}|{\rm sph}\,\DT\ein|^2
+2G_{_{\rm MX}}|{\rm dev}\,\DT\ein|^2
+\DT\alpha\,\eta'(\DT\alpha)
\\&
=\bbD_{_{\rm KV}}(\DT e{-}\DT\pi)\Colon(\DT e{-}\DT\pi)
+\bbD_{_{\rm MX}}\DT\pi\Colon\DT\pi+\DT\alpha\eta'(\DT\alpha)
\end{align}
when using also the notation (\ref{C-D}b,c).
Thus \eq{engr-abstract} reads more speficially here as 
\begin{align}\nonumber
&\frac{\d}{\d t}
\bigg(\int_\Omega\frac\varrho2|\DT u|^2+\frac32K_{_{\rm E}}|{\rm sph}(e(u){-}\pi)|^2
+
\calG_{_{\rm E}}(\alpha,{\rm dev}(e(u){-}\pi))
+\varrho\phi+\varrho\nabla\phi\Cdot u
\\\nonumber&\hspace{.2em}
+\frac\kappa2|\nabla\alpha|^2
+\frac12\varrho\big(|\omega{\cdot}(x{+}u)|^2{-}|\omega|^2|x{+}u|^2\big)
\,\d x
+\int_{\R^3}\frac1{8\uppi g}|\nabla\phi|^2\,\d x\bigg)
\\
&\hspace{.2em}
+\int_\Omega\!\bbD_{_{\rm KV}}(e(\DT u){-}\DT\pi)\Colon (e(\DT u){-}\DT\pi)
+\bbD_{_{\rm MX}}\DT\pi\Colon\DT\pi+\DT\alpha\partial\eta(\DT\alpha)\,\d x
=\int_{\R^3}\!\!\varrho_{\rm ext}(t)\DT\phi\,\d x.
\label{engr}\end{align}
Here one should naturally assume that $\eta$ is (naturally) nonsmooth possibly 
only at $\DT\alpha=0$ so that the dissipation rate 
$\DT\alpha\partial\eta(\DT\alpha)$ is actually 
a single-valued function.
The Coriolis force does not occur in \eq{engr}
because $\omega\times\DT u$ is always orthogonal to 
$\DT u$ so that $(\omega{\times}\DT u)\Cdot\DT u=0$. This is a well 
known effect that the Coriolis (pseudo) force does not make any work.

Denoting $\vec{n}$ the (unit) normal to the boundary $\partial\Omega$ of 
$\Omega$, i.e.\ the Earth surface, we complete the system by natural initial 
and boundary conditions: 
\begin{subequations}\label{BC-IC}\begin{align}
&u|_{t=0}=u_0,\ \ \ \ \DT u|_{t=0}=v_0,\ \ \ \ \pi|_{t=0}=\pi_0,\ \ \ \ \alpha|_{t=0}=\alpha_0.
\\
&\sigma\vec{n}=0\ \ \ \ \ \ \ \ \text{ on }\partial\Omega,
\ \ \ \ \phi(\infty)=0,\ \ \text{ and }\ 
\\&\label{BC-IC-alpha}
(\kappa\nabla\alpha){\cdot}\vec{n}=0\ \text{ on }\partial\Omega,
\end{align}\end{subequations}
Actually, it would make a good sense to consider also some period instead of initial conditions,
but the analysis of such a problem would be more difficult.

It is a conventional modelling property that $\alpha$ ranges the interval 
$[0,1]$, with $\alpha=0$ meaning no damage while $\alpha=1$ corresponds to a 
completely disintegrated rock. To make the damage model relatively simple 
without causing unnecessary analytic complications (leading e.g.\ to 
implementing the concept of nonsimple materials), the simplest modelling 
assumption that ensures the mentioned constraint $0\le\alpha\le1$ is 
$$
\partial_\alpha\calG_{_{\rm E}}(0,e)=0=\partial_\alpha\calG_{_{\rm E}}(1,e).
$$  
This facilitates the analysis of the model, sketched in Appendix below. Without 
going into details here, we can here summarize the main theoretical result as:

\begin{proposition}\label{prop-1}
The initial-boundary-value problem \eq{system+}--\eq{BC-IC} admits a weak 
solution in the sense of Definition~\ref{def-1} below.
This solution conserves energy in the sense that \eq{engr} holds when
integrated over any time interval $[0,t]$ with $0<t\le T$.
\end{proposition}

\begin{remark}[{\sl Activated creep: plasticity}]
\upshape
Some other modification of the ``solid'' model used in geophysical
literature modifies the evolution of the creep $\ein$ to be 
activated in the spirit what it standardly used in plasticity. Then
 $\zeta$ should be augmented by a term like $|{\rm dev}\ein|$ and,
to facilitate mathematical analysis, then $\varphi$ is to be augmented  
by the term $|\nabla{\rm dev}\ein|^2$. Then $\ein$ is called an inelastic
strain, rather than a Maxwellian creep. This model may be relevant in particular
to provide an additional dissipation of
energy important during big earthquakes and damage-dependent 
yield stress that may facilitate fast rupture.
\end{remark}

\begin{remark}[{\sl Geopotential}]
\upshape
In the standard presentation of the selfgraviting model,
the evolving potential $\phi$ is rather a perturbation of 
another gravitational potential (constant in time), resulted from 
a steady-state equilibrium configuration at a chosen initial time
with the Coriolis force 
zero and the centrifugal force 
$\varrho\omega{\times}(\omega{\times}x)$. The sum of these two is 
sometimes referred as a {\it geopotential}. Then another force occurs in 
\eq{system-u+}, causing the body $\Omega$ to be pre-stressed. Cf.\ e.g.\ 
\cite{Braz17EGEG,BrHoHo17VFEE,DahTro98TGS,WooDeu09TOEF}. 
\end{remark}


\section{Towards Boger viscoelastic fluid in the outer core and oceans}
\label{sec-Boger}

\begin{figure}
\psfrag{e}{\small $e=e(u)$}
\psfrag{s1}{\small $\sigma_{\rm sph}$}
\psfrag{s2}{\small $\sigma_{\rm dev}$}
\psfrag{e1}{\small ${\rm dev}\,e$}
\psfrag{e2}{\small ${\rm sph}\,e$}
\psfrag{G}{\small $G_{_{\rm E}}\!=G_{_{\rm E}}\!(\alpha,e)$}
\psfrag{Gkv}{\small $G_{_{\rm KV}}$}
\psfrag{K}{\small $K_{_{\rm E}}$}
\psfrag{Kkv}{\small $K_{_{\rm KV}}$}
\begin{center}
\includegraphics[width=19em]{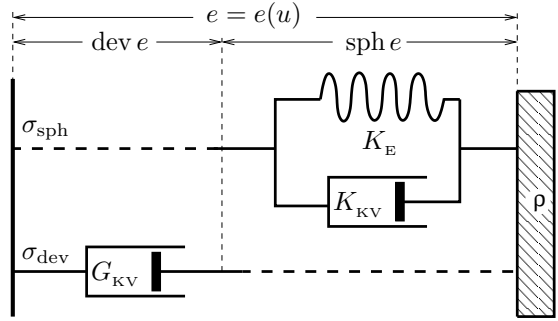}
\end{center}
\vspace*{-.6em}
\caption{\sl Schematic diagramme for the Boger-fluid rheology, resulting 
from Fig.~\ref{fig-mixed-response} when sending $K_{_{\rm MX}}\to\infty$, 
$G_{_{\rm MX}}\to\infty$, and $G_{_{\rm E}}\!(\alpha)\to0$. This model is considered 
in $\OF$ (=the outer core and oceans) coupled with the model from 
Figure~\ref{fig-mixed-response} in $\OS$ (=the mantle and the inner core).
}
\label{fig-mixed-response+}
\end{figure}

A standard specification of so-called Boger's viscoelastic fluids 
is as constant-viscosity elastic (non-Newtonian) fluids that 
behave as both liquids and solids. Originally, this concept was devised rather 
for materials like dilute polymer solutions. 
Yet, such (idealized) material may model the fluid outer core in the Earth,
which is a 2200\,km thick layer between the (rather solid) inner 
core of the radius about 1300\,km and the (rather solid) mantle, 
cf.~Fig.\,\ref{fig-earth}. The elasticity in the volumetric part is essential 
because it allows for a propagation of P-waves, which is a well documented 
phenomenon in the outer core, while S-waves are practically not penetrating 
this region because of its fluidic character as far as the shear response 
concerns, cf.\ Figure~\ref{fig-mixed-response+}.
Knowing quite reliably the speed $v_{_{\rm P}}$ of the P-waves 
and the mass density $\varrho$ depending on the depth, 
the elastic bulk modulus $K_{_{\rm E}}=K_{_{\rm E}}(x)$ in the outer core is easy to 
be seen from the formula $v_{_{\rm P}}=\sqrt{M/\varrho}$ where 
$M=\lambda+2G_{_{\rm E}}\!=K_{_{\rm E}}+\frac43G_{_{\rm E}}$ is the so-called P-wave modulus;
here $G_{_{\rm E}}=G_{_{\rm E}}(\alpha)$ in \eq{KV-limit+}. In the 
elastic fluid $G_{_{\rm E}}\!=0$ so that $M=K_{_{\rm E}}$. One can indeed determine simply as
$K_{_{\rm E}}=\varrho v_{_{\rm P}}^2$ when knowing $\varrho$ and $v_{_{\rm P}}$. 
According the reliably documented speed of the seismic P-waves
varying from 8\,km/s (top) to 10\,km/s (bottom) and the mass density 
varying from 10000\,kg/m$^3$ (top) to 12000\,kg/m$^3$ (bottom),
the elastic bulk modulus ranges from $K_{_{\rm E}}=640\,$GPa (top) to 
1200\,GPa (bottom).

For oceanic layers, it is relevant that the bulk modulus of seawater is about 
$K_{_{\rm E}}=M=2.3\,$GPa (increasing with pressure) which is, e.g., about 
80$\times$ smaller than in steel and about 5$\times$ smaller than in the crust, 
and thus water is remarkably elastically compressible. 
Modelling oceans as elastic fluid is thus relevant for propagation
seismic P-waves, which may manifest themselves on the surface by the 
Tsunami waves, cf.\ e.g.\ 
\cite{MeaFur13SSWO,TAKS13EEEE}. 
Note that, counting  the speed of P-waves close to the 
sea surface where $\varrho\sim10^3$\,kg/m$^3$ is 
$v_{_{\rm P}}=\sqrt{M/\varrho}\sim1.5\,$km/s, which is the sound speed,
so underwater acoustics exploits the same mechanism as seismic P-waves.

In addition to the elastic response, there is a (rather small) viscosity 
$(K_{_{\rm KV}},G_{_{\rm KV}})$ 
in the fluidic outer core and in oceans, too. In the outer core, some data 
indicates it of the order $10^{-2}$Pa\,s, cf.\ \cite{DKVD98VLIP}, while some 
uncertainty and big variation within depth $10^{4}$Pa\,s
seems also documented in literature, 
see e.g.\ \cite{Secc13VOC,SmyPal07VEOC}. The viscosity of the water
is even lower, of the order $10^{-3}$Pa\,s, varying with temperature and 
pressure. 

Here we can come to a viscoelastic (Boger's) fluid as a limit from the 
previous damageable Jeffreys' rheology when sending $G_{_{\rm E}}\!(\alpha)\to0$ 
so that, in particular, the resulting material naturally will 
not be subjected to any damage, and simultaneously 
$K_{_{\rm MX}}\to\infty$, and $G_{_{\rm MX}}\to\infty$.
As a result, $\pi\to0$ and eventually \eq{system-pi} is to hold only
on $\OS$ instead of the whole $\Omega$ while
\eq{system-alpha} still holds on the whole $\Omega$
but 
$\pl_\alpha\calG_{_{\rm E}}(\alpha,e_{\rm el})$
is zero on $\OF$
and values of the damage $\alpha$ are actually irrelevant 
 \eq{system-u+} on $\OF$.

More specifically, the system that results from \eq{system+}
by this limit procedure looks as:
\begin{subequations}\label{system++}\begin{align}\label{system-u++}
&\varrho\DDT u-{\rm div}(\sigma_{\rm sph}{+}\sigma_{\rm dev})
+\varrho(\omega{\times}(\omega{\times}(x{+}u))+2\omega{\times}\DT u+\nabla\phi)=0&&\text{in }\ \Omega
\\&\hspace{6.3em}\text{ with }\ 
\sigma_{\rm sph}=\begin{cases}3K_{_{\rm KV}}{\rm sph}\,
\DT e_{\rm el}
+3K_{_{\rm E}}{\rm sph}
\,e_{\rm el}\ \ \ \ \ \
&\text{in }\ \OS,\\
3K_{_{\rm KV}}{\rm sph}\,e(\DT u)
+3K_{_{\rm E}}{\rm sph}\,e(u)&\text{in }\ \OF,\end{cases}\hspace{-9em}
\\&\hspace{6.3em}\text{ and }\,\ \sigma_{\rm dev}=
\begin{cases}2G_{_{\rm KV}}{\rm dev}\,\DT e_{\rm el}
+
\partial_{e}\calG_{_{\rm E}}(\alpha,{\rm dev}\,e_{\rm el})
&\text{in }\ \OS,\\
2G_{_{\rm KV}}{\rm sph}\,e(\DT u)&\text{in }\ \OF,\end{cases}\hspace{-9em}
\\&\label{system-pi++}
\bbD_{_{\rm MX}}\DT \ein
=\sigma_{\rm sph}{+}\sigma_{\rm dev}
&&\text{in }\ \OS,
\\[.2em]&\label{system-alpha++}
\pl\eta(\DT\alpha)
+
\partial_\alpha G_{_{\rm E}}(\alpha,{\rm dev}\,e_{\rm el})
-{\rm div}(\kappa\nabla\alpha)\ni0
&&\text{in }\ \OS,
\\[-.3em]&\label{system4++}
{\rm div}\Big(\frac1{4\uppi g}\nabla\phi+\varrho u\Big)
=\varrho+\varrho_{\rm ext}(t)
&&\text{in }\ \R^3.
\end{align}\end{subequations}
The interface conditions for displacement/stress on the interior boundaries 
between the mantle and the 
outer core and the inner and the outer core as well between the 
mantle (crust) and the oceans are automatically involved as \eq{system-u++}
holds on the whole $\Omega$. In particular, the stress-vector equilibrium and 
continuity of the displacement across these interfaces are automatically
involved and does not need to be written explicitly. On the other hand,
note in particular that both $\pi$ and $\alpha$ in \eq{system++}
are now needed and defined only in the solid part $\OS$.
Therefore, the condition for the ``flux'' of $\alpha$ is now to be prescribed 
on all interior boundaries (i.e.\ on the mantle/core and mantle/oceans and
inner/outer-core interfaces). More specifically,
\eq{BC-IC-alpha} is to be replaced by
\begin{align}
\label{BC-IC-alpha+}
(\kappa\nabla\alpha){\cdot}\vec{n}=0\ \text{ on }\partial\OS.
\end{align}

We can justify this limit-passage scenario rigorously. To this goal, let us 
choose 
\begin{subequations}\label{G-K-eps}\begin{align}
&G_{_{\rm MX},\eps}(x)=\begin{cases}G_{_{\rm MX}}(x)/\eps,&\\\ G_{_{\rm MX}}(x),&\end{cases}
\ \ \;K_{_{\rm MX},\eps}(x)=\begin{cases}K_{_{\rm MX}}(x)/\eps\!\!&\text{for }x\In\OF,\\\ K_{_{\rm MX}}(x)\!\!&\text{for }x\In\OS,\end{cases}
\!\!\!
\\&\label{G-eps}
\calG_{_{\rm E},\eps}(x,\alpha,e)=\begin{cases}\eps\calG_{_{\rm E}}\!(x,\alpha,e),&\\\ \calG_{_{\rm E}}\!(x,\alpha,e),&\end{cases}\ 
\kappa_\eps(x)\  =\ \begin{cases}\eps\kappa(x)&\ \ \ \,\text{ for }x\In\OF,\\\ \kappa(x)&\ \ \ \,\text{ for }x\In\OS,\end{cases}
\\&\eta_\eps(x,\DT\alpha)\  =\ \begin{cases}\eps\eta(x,\DT\alpha)&
\ \ \ \text{ for }x\In\OF,
\\
\eta(x,\DT\alpha)&\ \ \ \text{ for }x\In\OS,\end{cases}
\end{align}\end{subequations}
while the other visco-elastic moduli $K_{_{\rm E}}$, $K_{_{\rm KV}}$, and 
$G_{_{\rm KV}}$ are kept fixed. The following statement will be made more specific
and proved in the Appendix:

\begin{proposition}\label{prop-2}
Let be some solution $(u_\eps,\ein_\eps,\alpha_\eps,\phi_\eps)$ 
of the initial-boundary-value problem \eq{system+}--\eq{BC-IC} 
with the data \eq{G-K-eps}, which do exist due to Proposition~\ref{prop-1}.
Then $(u_\eps,\ein_\eps|_{\QS}^{},\alpha_\eps|_{\QS}^{},\phi_\eps)$ 
converge (in terms of subsequences) for $\eps\to0$ to weak solutions to 
\eq{system++} with the initial/boundary conditions 
{\rm(\ref{BC-IC}a,b)}--\eq{BC-IC-alpha+}.
In particular, a weak solution to the initial-boundary-value problem
\eq{system++}--{\rm(\ref{BC-IC}a,b)}--\eq{BC-IC-alpha+} 
in the sense of Definition~\ref{def-2} below does exist.
In addition, this solution conserves energy.
Moreover, the energy dissipated through the Maxwell-type attenuation and 
by damage in
the fluidic regions $\OF$ over the time interval $I$ converges to zero, i.e.\ 
\begin{align}\label{Max-limit}
\!\int_{\QF}\!\!\frac32K_{_{\rm MX},\eps}|{\rm sph}\,\DT \ein_\eps|^2\!
+G_{_{\rm MX},\eps}|{\rm dev}\,\DT \ein_\eps|^2\!
+\DT\alpha_\eps\pl\eta(\DT\alpha_\eps)\,\d x\d t\to0\ \ \ \text{ for }\ 
\eps\to0.\
\end{align}
\end{proposition}

It should be emphasized that, although \eq{Max-limit} similarly as \eq{KV-limit}
below is intuitively quite clear and generally expected in geophysical 
literature,
its rigorous proof is not entirely trivial and relies on a regularity
which is rather automatic in linear systems but may be highly nontrivial or 
even false in nonlinear hyperbolic systems, as basically also here, 
cf.\ \cite{RajRou03EDSM,RoPaMa13QACV}.

\begin{remark}
\upshape
Instead of $G_{_{\rm MX}}\to\infty$ and $G_{_{\rm E}}\!(\alpha)\to0$,
one can think about sending $G_{_{\rm KV}}\to\infty$. Then, assuming 
the initial condition ${\rm dev}\,e_{\rm el}|_{t=0}=0$, we would have
 ${\rm dev}\,e_{\rm el}\to0$. Hence, in the limit we would see the model
from Fig.\,\ref{fig-mixed-response+} only with $G_{_{\rm KV}}$ replaced
by $G_{_{\rm MX}}$. Since the viscosity in fluidic regions is typically small
(as discuss in particular in the next section~\ref{sec4}) 
like also $G_{_{\rm KV}}$ while $G_{_{\rm MX}}$ is typically large in
solid-like geophysical materials, our choice \eq{G-K-eps} is
more straightforward.
Another difference would be that, instead of $\pi\to0$,
we would have only ${\rm sph}\,\pi\to0$ while 
${\rm dev}\,\pi\to{\rm dev}\,e(u)$.
\end{remark}

\begin{remark}
\upshape
Actually, the limit $K_{_{\rm MX}}\to\infty$ is relevant also in the solid
part $\OS$. This causes ${\rm sph}\,\pi\to0$ and in the limit one
obtains a model where $\pi$ is trace-free and accumulates only
shear strain, which is the most common ansatz in plasticity and creep, and in 
geophysical modelling, too.
\end{remark}

\begin{remark}[Incompressible Stokes-fluid limit]
\upshape
Passing simultaneously $K_{_{\rm E}}\to\infty$, we obtain in the limit the Stokes 
incompressible fluid. The bulk viscosity $K_{_{\rm KV}}$ then becomes irrelevant, 
so that only the shear viscosity $G_{_{\rm KV}}$ remains relevant in 
Figure~\ref{fig-mixed-response+}. Such model is a great idealization and
has limited application because, beside $v_{_{\rm S}}\to0$, now nonphysically 
$v_{_{\rm P}}\to\infty$. In particular, in this limit, ${\rm sph}\,e(u)$ is 
constant in time, namely ${\rm sph}\,e(u(x,t))=\frac13{\rm div}\,u_0(x)\bbI$.
Therefore ${\rm div}\,\DT u=0$, which expresses the usual incompressibility
condition. Such models are often used in geophysics in short-time-scale
models when the convective term $\varrho(\DT u\cdot\nabla)\DT u$, 
occurring in Navier-Stokes equation, can be neglected. 
Cf.\ e.g.\ \cite{ToCaMa09SSLV} for self-gravitating incompressible Stokes model
in layered geoid.
\end{remark}

\section{Suppressing viscosity towards purely elastic fluid}\label{sec4}

As already said, the viscosity $(K_{_{\rm KV}},G_{_{\rm KV}})$ 
in the fluidic domains is only very small (and even not certainly known 
in deep parts of the outer core). Note also that, when $G_{_{\rm E}}\!=0$
in the fluidic domains, $G_{_{\rm KV}}$ is in the position of the 
Maxwellian viscosity. For example, in mantle, $10^2-10^4\,$Pa\,s is 
considered in \cite{LyHaBZ11NLVE}.
As for the fluidic regions, as said above, the oceans exhibit viscosity around 
$10^{-3}$Pa\,s while the outer core around $10^{-2}-10^{4}$Pa\,s.
In any case, these viscosities are much smaller than
typical viscosity in the crust which is of the order $10^{22}-10^{24}\,$Pa\,s
or in the inner core of the order about $10^{14}-10^{15}\,$Pa\,s, 
see e.g.\ \cite{KooDum11VEIC}.
This certainly gives a good motivation to study the asymptotics when this 
viscosity goes to zero.
In the limit, it yields an inviscid, purely elastic model where the 
Hooke elasticity counts only with the spherical response while the 
shear-stress-free response imitates the ideal inviscid fluid;
sometimes, such sort of models are called compressible Euler fluids.
Such materials are called {\it elastic fluids},
or more specifically just compressible inviscid fluids.
These fluids allow still for propagation of 
{\it P-waves}\index{P-wave!in elastic fluids} 
(whose speed is $\varrho^{-1/2}K_{_{\rm E}}^{1/2}$ as already mentioned) 
while S-waves are completely excluded.

The resulted system is again (\ref{system++}a,d-f) but now with 
\begin{subequations}\label{system+++}\begin{align}
&\quad
\sigma_{\rm sph}=\begin{cases}3K_{_{\rm KV}}{\rm sph}(e(\DT u){-}\DT\pi)
+3K_{_{\rm E}}{\rm sph}(e(u){-}\pi)\ \ \ \ \ \ 
&\text{in }\OS,\\
3K_{_{\rm E}}{\rm sph}\,e(u)&\text{in }\OF,\end{cases}
\\&\quad
\sigma_{\rm dev}=
\begin{cases}2G_{_{\rm KV}}{\rm dev}(e(\DT u){-}\DT\pi)
+
\partial_e^{}\calG_{_{\rm E}}\!(\alpha,{\rm dev}(e(u){-}\pi))
&\text{in }\OS,\\
0&\text{in }\OF.\end{cases}
\end{align}\end{subequations}

The weak formulation for the limit problem as far as the 
force equilibrium 
(\ref{system++}a,d) 
with \eq{system+++} 
concerns 
arises when testing \eq{system-u++} by a smooth test function 
$v$ with $v(T)=0$ and ${\rm dev}\,e(v)=0$ on $\QF$ and making 
by part integration in time and applying Green formula in space, cf.\
Definition~\ref{def-3} in the Appendix.

We can justify this limit-passage scenario rigorously. To this goal, let us 
choose 
\begin{align}\label{KV-eps}
K_{_{\rm KV},\eps}(x)=\begin{cases}\eps K_{_{\rm KV}}(x),\!\!&
\\\ K_{_{\rm KV}}(x),\!\!&\end{cases}\ \ \ 
G_{_{\rm KV},\eps}(x)=\begin{cases}\eps G_{_{\rm KV}}(x)\!\!&\text{for }x\In\OF,
\\\ G_{_{\rm KV}}(x)\!\!&\text{for }x\In\OS,\end{cases}
\end{align}
while $G_{_{\rm E}}\!(\alpha)=0$ and $K_{_{\rm MX}}=G_{_{\rm MX}}=\infty$ on $\OF$, resulting
already from the limit in Section~\ref{sec-Boger}.
The following statement
will be proved in the Appendix:

\begin{proposition}\label{prop-3}
If the initial conditions are enough smooth (cf.\ the assumption in
Definition~\ref{def-3}
below),
the solutions $(u_\eps,\ein_\eps,\alpha_\eps,\phi_\eps)$ 
of the initial-boundary-value problem 
\eq{system++}--{\rm(\ref{BC-IC}a,b)}--\eq{BC-IC-alpha+}
with the data \eq{KV-eps}, which does exist due to Proposition~\ref{prop-2},
converge (in terms of subsequences) for $\eps\to0$ to weak solutions to 
{\rm(\ref{system++}a,d-f)}--{\rm(\ref{BC-IC}a,b)}--\eq{BC-IC-alpha+}
with \eq{system+++}.
In particular, a weak solution to this initial-boundary-value problem
in the sense of Definition~\ref{def-3} below does exist. Moreover,
this solution conserves energy
and the energy dissipated through the Kelvin-Voigt-type attenuation in
the fluidic regions $\OF$ over the time interval $I$ converges to zero, i.e.\ 
\begin{align}\label{KV-limit}
\eps\int_{\QF}\!\!\frac32K_{_{\rm KV}}|{\rm sph}\,e(\DT u_\eps)|^2\!
+G_{_{\rm KV}}|{\rm dev}\,e(\DT u_\eps)|^2\,\d x\d t\to0
\ \ \ \text{ for }\ 
\eps\to0.\
\end{align}
\end{proposition}

\section{Appendix: analysis sketched}\label{sec-anal}

In what follows, we will use the (standard) notation for the Lebesgue 
$L^p$-spaces and $W^{k,p}$ for Sobolev spaces whose $k$-th distributional 
derivatives are in $L^p$-spaces, and the abbreviation $H^k=W^{k,2}$. 
Moreover, we will use the standard notation $p'=p/(p{-}1)$.
In the vectorial case, we will write $L^p(\Omega;\R^3)\cong L^p(\Omega)^3$ 
and $W^{1,p}(\Omega;\R^3)\cong W^{1,p}(\Omega)^3$. 
%
For the fixed time interval $I=[0,T]$, we denote by $L^p(I;X)$ the 
standard Bochner space of Bochner-measurable mappings $I\to X$ with $X$ a 
Banach space. Also, $W^{k,p}(I;X)$ denotes the Banach space of mappings 
from $L^p(I;X)$ whose $k$-th distributional derivative in time is also in 
$L^p(I;X)$. The dual space to $X$ will be denoted by $X^*$.
The scalar product between vectors, matrices, or 3rd-order tensors
will be denoted by ``$\,\cdot\,$'', ``$\,:\,$'', or ``$\,\Vdots\,$'', 
respectively. Finally, in what follows, $K$ denotes a positive, possibly 
large constant.


We will impose the basic assumptions on the data and incorporate them
directly into Definitions~\ref{def-1}-\ref{def-3}. In particular,
\begin{subequations}\label{ass}\begin{align}\label{ass-rho}
&{\rm inf}_{x\in\Omega}^{}\varrho(x)>0,\ \ \ \ \varrho=0\ \text{ on }\R^3{\setminus}\Omega,\ \ \ \ \varrho\in L^\infty(\Omega)\cap W^{1,3}(\Omega),
\\&\varrho_{\rm ext}\In W^{1,1}(I;L^\infty(\R^3)),\ \ 
\varrho_{\rm ext}=0
\text{ outside a bounded set in }I{\times}\R^3,
\\&\nonumber
K_{_{\rm E}},G_{_{\rm KV}},K_{_{\rm KV}},G_{_{\rm MX}},K_{_{\rm MX}},\calG_{_{\rm E}}(\cdot,\alpha,e),
\kappa:\Omega\to[0,\infty)\ \text{ measurable and}
\\&\nonumber\qquad\qquad
\!\inf_{x\in\Omega}
\min\big(G_{_{\rm KV}}(x),K_{_{\rm KV}}(x),G_{_{\rm MX}}(x),K_{_{\rm MX}}(x),K_{_{\rm E}}(x),\kappa\big)
>0,\ \ \ 
\\&\qquad\qquad\label{ass-small-KV/MX}
\sup_{x\in\OS}\max\big(K_{_{\rm KV}}(x)/K_{_{\rm MX}}(x),G_{_{\rm KV}}(x)/G_{_{\rm MX}}(x)\big)<1,
\\&\nonumber
\calG_{_{\rm E}}(x,\cdot,\cdot){:}
[0,1]{\times}\R^{3\times3}\to[0,\infty)\ 
\text{ twice continuously differentiable},\ \
\\&\nonumber\qquad\qquad
\calG_{_{\rm E}}(x,\alpha,\cdot){:}\R^{3\times3}\to[0,\infty)\ \text{ convex},
\\&\nonumber\qquad\qquad
|\pl_e\calG_{_{\rm E}}(x,\alpha,e)|\le C(1+|e|),\ \ \ |\pl_\alpha\calG_{_{\rm E}}(x,\alpha,e)|\le C(1+|e|),\ \ \ 
\\&\label{ass-growth-G}\qquad\qquad
|\pl_{e e}^2\calG_{_{\rm E}}(x,\alpha,e)|\le C,\ \ \ \ \ \
|\pl_{\alpha e}^2\calG_{_{\rm E}}(x,\alpha,e)|\le C,
\\&\eta:\R\to\R\ \ \text{ uniformly convex},\ \eta(0)=0,
\\\nonumber&u_0\In H^1(\Omega;\R^3),\ \ 
v_0\In L^2(\Omega;\R^3),\ \ \pi_0\In L^2(\Omega;\R^{3\times 3}),\ \ 
\pi_0=0\ \text{ on }\,\OF,
\\&\qquad\qquad\qquad\text{ and }
\alpha_0\In H^1(\Omega)\ 
\text{ with $0\le\alpha_0(x)\le1$ for all $x\In\Omega$.
}
\end{align}\end{subequations}
Let us emphasize that \eq{ass-small-KV/MX} is
used only for the second limit passage for estimation 
\eq{est-regular}--\eq{est-regular+} and
is well satisfied in geophysical models where the ratio
of the Kelvin-Voigt and Maxwell viscosities in solid regions is surely
below $10^{-8}$, 
as said in Section~\ref{sec4}. Also, let us emphasize that
the growth restrictions imposed on 
$\pl_\alpha^{}\calG_{_{\rm E}}$ and 
$\pl_{\alpha e}^2\calG_{_{\rm E}}$ in \eq{ass-growth-G}
are compatible with the ansatz \eq{KV-limit+} and used 
in \eq{system-alpha-w} and for \eq{est-regular} below.

We further
integrate \eq{system-alpha-w-} over $Q$ and apply the Green theorem and the by-part
integration in time to the term $(v{-}\DT\alpha){\rm div}(\kappa\nabla\alpha)$,
cf.\ Remark~\ref{rem-flow-rule}. 

\begin{definition}[Weak solution to \eq{system+}--\eq{BC-IC}]\label{def-1}
The quadruple $(u,\pi,\alpha,\phi)$ is called a weak solution to 
the initial-boundary-value problem \eq{system+}--\eq{BC-IC} 
provided the data satisfies \eq{ass} and 
\begin{subequations}\label{est}\begin{align}\label{est-u}
&u\In W^{1,\infty}(I;L^2(\Omega;\R^3))\cap H^1(I;H^1(\Omega;\R^3)),
\\&\pi\In H^1(I;L^2(\Omega;\R^{3\times 3})),
\\&\label{est-a}
\alpha\in L^\infty(I;H^1(\Omega))\cap H^1(I;L^2(\Omega)),
\\&\label{est-phi}
\phi\In L^\infty(I;H^1(\R^3)),
\end{align}\end{subequations} 
and the integral identity 
\begin{subequations}\label{system-w}\begin{align}\label{momentum-w}
&\int_Q\sigma{:}e(v)-\varrho\DT u{\cdot}\DT v+f{\cdot}v\,\d x\d t=\int_\Omega
\varrho v_0{\cdot}v(0)\,\d x
\intertext{holds for any $v\in H^1(Q;\R^3)$ with $v(T)=0$ and with $\sigma$ and $f$ from \eq{system-u+},}
&\nonumber\int_Q\big(v{-}\DT\alpha\big)
\partial_\alpha\calG_{_{\rm E}}(\alpha,{\rm dev}(e(u){-}\pi))
+\kappa\nabla\alpha\Cdot\nabla v+\eta(v)\,\d x\d t
\\[-.4em]&\qquad\qquad\qquad
+\int_\Omega\frac\kappa2|\nabla\alpha_0|^2\,\d x
\ge\int_Q\!\eta(\DT\alpha)\,\d x\d t
+\int_\Omega\frac\kappa2|\nabla\alpha(T)|^2\,\d x
\label{system-alpha-w}
\intertext{for all $v\in L^2(I;H^1(\Omega))$, 
}
&\int_{I\times\R^3}
\Big(\frac{1}{4\pi g}\nabla\phi+\varrho u\Big){\cdot}\nabla v+
\big(\varrho+\varrho_{\rm ext}{}\big)v\,\d x\d t=0
\label{gravity-w}
\end{align}\end{subequations} 
for all $v\in L^2(I;H^1(\R^3))$, 
and also \eq{system-pi} holds a.e.\ on $Q$, and eventually the resting 
initial conditions $u(0)=u_0$, $\pi(0)=\pi_0$, and $\alpha(0)=\alpha_0$ hold 
a.e.\ on $\Omega$.
\end{definition}

\noindent{\it Sketched proof of Proposition~\ref{prop-1}.}
Without going into details, we may expect that we applied some approximation
method (e.g.\ a Galerkin-type approximation in space) to obtain some
approximate solution to the 
initial-boundary-value problem \eq{system+}--\eq{BC-IC}
which, in principle, can be implemented on 
computers e.g.\ by the finite-element method. Then the a-priori 
estimates in the space as in \eq{est} hold also for the approximate
solutions, which may be interpreted as a numerical stability 
of the specific approximation scheme. This approximation leads (after 
smoothening of the potential $\eta$) to the initial-value problem for a
system of ordinary-differential equations. Existence of its solution, let
us denote it by $(u_k,\pi_k,\alpha_k,\phi_k)$ with $k\in\N$ referring to the 
finite-dimensional subspaces used for the Galerkin discretisation, can be 
proved by successive-continuation argument, using the $L^\infty(I)$-estimates 
below. Here we may assume that the initial conditions lie in the 
finite-dimensional spaces used for the Galerkin approximation so that no
further approximation is needed. 

The energy balance \eq{engr} can serve to see basic apriori estimates and
to perform analysis of the system \eq{system+}.
Integrating \eq{engr} over $[0,t]$ and using the by-part integration in time 
for the power of tidal load $\varrho_{\rm ext}\DT\phi$ and the notation 
(\ref{C-D}b,c), we write
\begin{align}\nonumber
&
\int_\Omega\bigg(\frac\varrho2|\DT u_k(t)|^2+
\frac32K_{_{\rm E}}|{\rm sph}
\,e_{{\rm el},k}(t)|^2
+\calG_{_{\rm E}}(\alpha_k(t),{\rm dev}\,e_{{\rm el},k}(t))
\\\nonumber&\hspace{2em}
+\frac\kappa2|\nabla\alpha_k(t)|^2
+
\frac12\varrho|\omega{\cdot}(x{+}u_k(t))|^2\bigg)\,\d x
+\int_{\R^3}\frac1{
8\uppi g}|\nabla\phi_k(t)|^2\,\d x
\\\nonumber
&\hspace{3em}
+\int_0^t\!\!\int_\Omega\bbD_{_{\rm KV}}
\,\DT e_{{\rm el},k}\Colon\DT e_{{\rm el},k}
+\bbD_{_{\rm MX}}\DT\pi_k\Colon\DT\pi_k+\DT\alpha_k\eta'(\DT\alpha_k)\,\d x\d t
\\&\nonumber\hspace{7em}
\le\int_\Omega
-\varrho\nabla\phi_k(t)\Cdot u_k(t)-\varrho\phi_k(t)\,\d x
\\[-.4em]&\hspace{9em}+\int_{\R^3}\varrho_{\rm ext}(t)\phi_k(t)\,\d x
-\int_0^t\!\!\int_{\R^3}\DT\varrho_{\rm ext}\phi_k\,\d x\d t+E_0
\label{engr-int}\end{align}
with $e_{{\rm el},k}=e(u_k)-\pi_k$ and with 
the upper bound for the initial energy 
\begin{align}\nonumber E_0&=\int_\Omega
\frac\varrho2|v_0|^2+
\frac32K_{_{\rm E}}|{\rm sph}
\,e_{{\rm el},0}|^2
+\calG_{_{\rm E}}(\alpha_0,e_{{\rm el},0})
+\varrho\phi(0)+\varrho\nabla\phi(0)\Cdot u_0
\\[-.4em]
&\qquad\ \ \ \
+\frac12\varrho|\omega{\cdot}(x{+}u_0)|^2
+\frac\kappa2|\nabla\alpha_0|^2\,\d x
+\int_{\R^3}\!\varrho_{\rm ext}(0,x)\phi(0,x)\,\d x,
\label{def-E0}
\end{align}
where $e_{{\rm el},0}=e(u_0)-\pi_0$ and where
$\phi(0,x)\in H^1(\Omega)$ is the gravitational potential 
solving \eq{system4} for $t=0$. The inequality in \eq{engr-int} have arisen 
from the energy equality \eq{engr} written for the approximate solution 
and integrated over the time interval $[0,t]$ by forgetting 
some parts of the centrifugal potential with a guaranteed sign.
The non-coercive contribution of the centrifugal force on
the right-hand side of \eq{engr-int} is to be estimated by H\"older's inequality
as 
\begin{align}\nonumber
\int_\Omega
\frac12\varrho|\omega|^2|x{+}u_k(t)|^2
\,\d x
&=\int_\Omega\frac12\varrho|\omega|^2
\bigg|x{+}u_0(x){+}\int_0^t\DT u_k(\tau,x)\,\d\tau\bigg|^2\,\d x
\\&
\le C\bigg(1+\int_0^t\|\DT u_k(\tau)\|_{L^2(\Omega;\R^3)}^2\,\d\tau\bigg)
\end{align}
with some constant $C$, and then treated by Gronwall's inequality, exploiting 
the kinetic energy on the left-hand side of \eq{engr-int}.
Furthermore, we can estimate $\int_\Omega\nabla\phi_k{\cdot}u_k\,\d x
\le\epsilon\|\nabla\phi_k\|_{L^2(\Omega;\R^3)}^2
+C_\epsilon(\int_0^t\|\DT u_k\|_{L^2(\Omega;\R^3)}^2\,\d t+\|u_0\|_{L^2(\Omega;\R^3)}^2)$.
Also use the estimate $\|\phi_k\|_{L^2(\Omega)}\le C\|\nabla\phi_k\|_{L^2(\R^3;\R^3)}$ 
relying on the boundedness of $\Omega$, provided $\phi(\infty)=0$  
which is a standard ``boundary''
condition for the gravitational potential used in geophysics.
The last integral on the right-hand side of \eq{engr-int} bears the estimation
\begin{align}\nonumber
&\int_0^t\!\!\int_{\R^3}\DT\varrho_{\rm ext}\phi_k\,\d x\d t\le
\int_0^t\!\|\DT\varrho_{\rm ext}\|_{L^\infty(\R^3)}\|\phi_k\|_{L^1(\R^3)}\d t
\\[-.3em]&\le C\int_0^t\!\|\DT\varrho_{\rm ext}\|_{L^\infty(\R^3)}\|\nabla\phi_k\|_{L^2(\R^3;\R^3)}\d t
\le\int_0^t\!\|\DT\varrho_{\rm ext}\|_{L^\infty(\R^3)}\big(1+\|\phi_k\|_{H^1(\R^3)}^2\big)\,\d t,
\nonumber\end{align}
where we used $\|\phi_k\|_{
L^1(\R^3)}\le C\|\nabla\phi_k\|_{L^2(\R^3;\R^3)}$;
here the ``boundary'' condition $\phi_k(\infty)=0$ together with 
$\varrho+\varrho_{\rm ext}$ compactly supported is used.
Altogether, using the Gronwall inequality, we obtain the a-priori estimates 
for the approximate solution:
\begin{subequations}\label{est+}\begin{align}\label{est-u+}
&\|u_k\|_{W^{1,\infty}(I;L^2(\Omega;\R^3))\cap 
H^1(I;H^1(\Omega;\R^3))}^{}\le C,
\\&\label{est-pi}
\|\pi_k\|_{H^1(I;L^2(\Omega;\R^{3\times 3}))}^{}\le C,
\\&\label{est-a+}
\|\alpha_k\|_{L^\infty(I;H^1(\Omega))\cap H^1(I;L^2(\Omega))}^{}\le C,
\\&\|\phi_k\|_{L^\infty(I;H^1(\R^3))}^{}\le C.
\end{align}\end{subequations} 
By Banach's selection principle, we consider a weakly* convergent subsequence
respecting the topologies specified in \eq{est+}. 
For the limit passage in the nonlinear term $\pl_\alpha\calG_{_{\rm E}}$ in 
\eq{system-alpha-w}, we need to improve it for the 
strong convergence of $e_{{\rm el},k}\to e_{\rm el}$ in $L^2(Q;\R^{3\times3})$. 
To this goal, we use the test function $u_k{-}u$ and $\pi_k{-}\pi$ for 
the Galerkin approximation of \eq{system-u+} and \eq{system-pi}, respectively,
integrate it over the time interval $[0,t]$, and estimate
\begin{align}\nonumber
&\!\int_\Omega\frac12\bbD_{_{\rm KV}}(e_{{\rm el},k}(t){-}e_{\rm el}(t))\Colon(e_{{\rm el},k}(t){-}e_{\rm el}(t))
\d x
\\[-.4em]\nonumber
&\!\le\frac12\int_\Omega\!\!\bbD_{_{\rm KV}}(e_{{\rm el},k}(t){-}e_{\rm el}(t))\Colon(e_{{\rm el},k}(t){-}e_{\rm el}(t))
{+}\bbD_{_{\rm MX}}(\pi_k(t){-}\pi(t))\Colon(\pi_k(t){-}\pi(t))
\d x
\\[-.1em]&\nonumber
+\int_0^t\!\!\int_\Omega\!\!K_{_{\rm E}}|{\rm sph}(e_{{\rm el},k}{-}e_{\rm el})|^2
+\big(\pl_e^{}\calG_{_{\rm E}}(\alpha_k,e_{{\rm el},k})
-\pl_e^{}\calG_{_{\rm E}}(\alpha_k,e_{\rm el})\big)
\Colon(e_{{\rm el},k}{-}e_{\rm el})\,\d x\d t
\\&\nonumber=\int_0^t\!\!\int_\Omega-(\varrho\DDT u_k+f_k)\Cdot(u_k{-}u)
\\&\nonumber\qquad-\big(\bbD_{_{\rm KV}}\DT e_{\rm el}+K_{_{\rm E}}{\rm sph}e_{\rm el}
+\pl_e^{}\calG_{_{\rm E}}(\alpha_k,e_{\rm el})\big)\Colon(e_{{\rm el},k}{-}e_{\rm el})
-\bbD_{_{\rm MX}}\DT\pi\Colon(\pi_k{-}\pi)\,\d x\d t
\displaybreak\\&\nonumber=
\int_0^t\!\!\int_\Omega \varrho\DT u_k\Cdot(\DT u_k{-}\DT u)-f_k\Cdot(u_k{-}u)
+\big(\pl_e^{}\calG_{_{\rm E}}(\alpha,e_{\rm el}){-}\pl_e^{}\calG_{_{\rm E}}(\alpha_k,e_{\rm el})
\big)\Colon(e_{{\rm el},k}{-}e_{\rm el})
\\&\nonumber\qquad\ -\big(\bbD_{_{\rm KV}}\DT e_{\rm el}+K_{_{\rm E}}{\rm sph}e_{\rm el}
+\pl_e^{}\calG_{_{\rm E}}(\alpha,e_{\rm el})\big)\Colon(e_{{\rm el},k}{-}e_{\rm el})
-\bbD_{_{\rm MX}}\DT\pi\Colon(\pi_k{-}\pi)\,\d x\d t
\\[-.2em]&\qquad\qquad\qquad\qquad\qquad\qquad\qquad\qquad
-\int_\Omega\varrho\DT u_k(t)\Cdot(u_k(t){-}u(t))\,\d x\ \to\ 0\,
\label{KV-damage-strong-e(t)}
\end{align}
with $f_k=\varrho(\omega{\times}(\omega{\times}(x{+}u_k))+
2\omega{\times}\DT u_k+\nabla\phi_k)$; note that $f_k$ is 
bounded in $L^2(Q;\R^3)$ due to the a-priori estimates (\ref{est+}a,d).
Actually, \eq{KV-damage-strong-e(t)} is to be understood rather as 
a conceptual strategy:
the mentioned test functions $u_k{-}u$ and $\pi_k{-}\pi$
are not legitimate for the Galerkin approximation and $(u,\pi)$ is still
to be approximated strongly to be valued in the respective finite-dimensional
subspaces - we omitted these standard technical details for simplicity. 
For the convergence to 0 in \eq{KV-damage-strong-e(t)}, we used that 
$\DT u_k\to\DT u$ strongly $L^2(Q;\R^3)$ due to
the Aubin-Lions theorem, relying on the estimate \eq{est-u+} together
with an information about $\DDT u$ from the equation \eq{system-u+} itself,
and also $\DT u_k(t)$ is bounded in $L^2(\Omega;\R^3)$
while  $u_k(t)\to u(t)$ strongly in $L^2(\Omega;\R^3)$
due to the Rellich compact embedding $H^1(\Omega)\subset L^2(\Omega)$,
and also 
$\pl_e^{}\calG_{_{\rm E}}(\alpha_k,e_{\rm el})\to\pl_e^{}\calG_{_{\rm E}}(\alpha,e_{\rm el})$
strongly in $L^2(Q;\Rdev)$ 
due to the continuity of the Nemytski\u\i\ mapping 
induced by $\pl_e^{}\calG_{_{\rm E}}(\cdot,e_{\rm el})$ and 
$\alpha_k\to\alpha$ in $L^2(Q)$ again just by the Rellich theorem
since both $\DT\alpha_k$ and $\nabla\alpha_k$ is estimated in $L^2$-spaces.
Thus, from \eq{KV-damage-strong-e(t)}, we obtain
$e_{{\rm el},k}(t)\to e_{\rm el}(t)$ strongly in $L^2(\Omega;\R^{3\times 3})$ 
for all $t\in I$. Using 
it for a general $t\In I$, we obtain $e_{{\rm el},k}\to e_{\rm el}$ strongly in 
$L^2(Q;\R^{3\times 3})$ by the Lebesgue theorem.


The limit passage in the Galerkin approximation towards the 
integral identities \eq{system-w} is then simple by weak/strong
continuity or, in case of the variational inequality \eq{system-alpha-w}, 
also semicontinuity.

For the energy conservation, the essential needed facts are that 
$\sqrt\varrho\DDT u\in L^2(I;H^1(\Omega;\R^3)^*)$ is in duality with 
$\sqrt\varrho\DT u\in L^2(I;H^1(\Omega;\R^3))$ and also that 
${\rm div}(\kappa\nabla\alpha)\in L^2(Q)$ is in duality with 
$\DT\alpha\in L^2(Q)$ so that the by-part
integration formulas rigorously hold:
\begin{subequations}\begin{align}\label{calculus1}
&\int_0^t\big\langle\sqrt\varrho\DDT u,\sqrt\varrho\DT u\big\rangle\,\d t
=\int_\Omega\frac\varrho2|\DT u(t)|^2-\frac\varrho2|\DT u(0)|^2\,\d x,\ \ \ 
\text{ and}
\\&\label{calculus2}
\int_0^t\!\!\int_\Omega{\rm div}(\kappa\nabla\alpha)\DT\alpha\,\d x\d t
=\int_\Omega\frac\kappa2|\nabla\alpha(0)|^2
-\frac\kappa2|\nabla\alpha(t)|^2\,\d x\,;
\end{align}\end{subequations} 
see e.g.\ \cite[Formulas (7.22) and (12.133b)]{Roub13NPDE}. More in detail, for 
\eq{calculus1} we have used the comparison 
$\sqrt\varrho\DDT u=({\rm div}\,\sigma-f)/\sqrt\varrho$
and the estimate 
\begin{align}\nonumber
&\big\|\sqrt\varrho\DDT u\big\|_{L^2(I;H^1(\Omega;\R^3)^*)}
=
\sup_{\|v\|_{L^2(I;H^1(\Omega;\R^3))}\le1}
\int_Q\frac{{\rm div}\,\sigma{-}f}{\sqrt\varrho}v\,\d x\d t
\\&\qquad\qquad
=\sup_{\|v\|_{L^2(I;H^1(\Omega;\R^3))}\le1}\int_Q\sigma\Colon\Big(\sqrt\varrho\nabla v-
\frac{v{\otimes}\nabla\varrho}{\sqrt\varrho}\Big)
+\frac{f}{\sqrt\varrho}v\,\d x\d t
\le C
\label{est-of-DDTu}\end{align}
for which the smoothness $\varrho\in W^{1,3}(\Omega)$ is needed due to
the occurrence of $\nabla\varrho$ in \eq{est-of-DDTu}, 
while for \eq{calculus2}
we have used the comparison ${\rm div}(\kappa\nabla\alpha)\in\pl\zeta(\DT\alpha)
+\pl_\alpha\calG_{_{\rm E}}(\alpha,{\rm dev}\,e_{\rm el})$.
$\Box$

\begin{definition}[Weak solution to 
{\rm(\ref{BC-IC}a,b)}--\eq{system++}--\eq{BC-IC-alpha+}.]\label{def-2}
The quadruple $(u,\pi,\alpha,\phi)$ is called a weak solution to 
the initial-boundary-value problem  
{\rm(\ref{BC-IC}a,b)}--\eq{system++}--\eq{BC-IC-alpha+}
with the data \eq{G-K-eps} provided again the data satisfies \eq{ass} and 
(\ref{est+}a,d) holds together with 
\begin{align}\label{est-2-}
\pi\In H^1(I;L^2(\OS;\Rsym))
\ \ \text{ and }\ \ 
\alpha\in L^\infty(I;H^1(\OS))\cap H^1(I;L^2(\OS))
\end{align}
and the integral identity \eq{momentum-w} holds with 
$\sigma=\sigma_{\rm sph}+\sigma_{\rm dev}$ from (\ref{system++}b,c), 
also \eq{system-alpha-w} holds with 
$Q$ and $\Omega$ replaced respectively $\QS$ and $\OS$,
furthermore \eq{gravity-w} holds, 
and also \eq{system-pi} holds a.e.\ on $\QS$, and eventually 
$u(0)=u_0$ holds a.e.\ on $\Omega$, and $\pi(0)=\pi_0$ and 
$\alpha(0)=\alpha_0$ hold a.e.\ on $\OS$.
\end{definition}

\noindent{\it Sketched proof of Proposition~\ref{prop-2}.}
Like \eq{engr-int} but now with \eq{G-K-eps} taken into account, we have
\begin{align}\nonumber
&
\int_\Omega\frac\varrho2|\DT u_\eps(t)|^2
+\frac32K_{_{\rm E}}|{\rm sph}\,e_{{\rm el},\eps}(t)|^2+
\frac12\varrho|\omega{\cdot}(x{+}u_\eps(t))|^2\,\d x
+\int_{\R^3}\frac1{
8\uppi g}|\nabla\phi_\eps(t)|^2\,\d x
\\[-.4em]\nonumber&\hspace{2em}
+\int_{\OS}\calG_{_{\rm E}}(\alpha_\eps(t),{\rm dev}\,e_{{\rm el},\eps}(t))
+\frac\kappa2|\nabla\alpha_\eps(t)|^2\,\d x
\\[-.4em]\nonumber&\hspace{9em}+
\eps\int_{\OF}\calG_{_{\rm E}}(\alpha_\eps(t),{\rm dev}\,e_{{\rm el},\eps}(t))
+\frac\kappa2|\nabla\alpha_\eps(t)|^2\,\d x
\\\nonumber
&\hspace{2em}
+\int_0^t\!\!\int_\Omega\bbD_{_{\rm KV}}
\,\DT e_{{\rm el},\eps}\Colon\DT e_{{\rm el},\eps}\,\d x\d t
+\int_0^t\!\!\int_{\OS}
\bbD_{_{\rm MX}}\DT\pi_\eps\Colon\DT\pi_\eps
+\DT\alpha_\eps\eta'(\DT\alpha_\eps)\,\d x\d t
\\[-.4em]\nonumber&\hspace{2em}
+\int_0^t\!\!\int_{\OF}
\frac1\eps\bbD_{_{\rm MX}}\DT\pi_\eps\Colon\DT\pi_\eps
+\eps\DT\alpha_\eps\eta'(\DT\alpha_\eps)\,\d x\d t
\le\int_\Omega
-\varrho\nabla\phi_\eps(t)\Cdot u_\eps(t)-\varrho\phi_\eps(t)\,\d x
\\[-.4em]&\hspace{9em}+\int_{\R^3}\varrho_{\rm ext}(t)\phi_\eps(t)\,\d x
-\int_0^t\!\!\int_{\R^3}\DT\varrho_{\rm ext}\phi_\eps\,\d x\d t+E_0
\label{engr-int-2}\end{align}
with $E_0$ again from \eq{def-E0}, relying that $\eps\le1$.
By the Gronwall-inequality arguments like used for \eq{engr-int},
we obtain the a-priori estimates (\ref{est+}b,d) now for $(\pi_\eps,\phi_\eps)$
instead of $(\pi_k,\phi_k)$ and also 
\begin{subequations}\label{est-2}\begin{align}\label{est-u-2}
&\|u_\eps\|_{W^{1,\infty}(I;L^2(\Omega;\R^3))\cap 
H^1(I;H^1(\OS;\R^3))}^{}\le C
\\&\label{est-a-2}
\|\alpha_\eps\|_{L^\infty(I;H^1(\OS))\cap H^1(I;L^2(\OS))}^{}\le C,
\\&\label{est-pi-2}
\|\DT\pi_\eps\|_{L^2(\QF;\R^{3\times 3})}^{}\le\sqrt\eps C,
\\&\label{est-DT-a-2}
\|\DT\alpha_\eps\|_{L^2(\QF;\R^{3\times 3})}^{}\le C/\sqrt\eps,
\ \ \ \text{ and }\ \ \ 
\\&\label{est-nabla-a-2}
\|\nabla\alpha_\eps\|_{L^\infty(I;L^2(\OF;\R^3))}\le C/\sqrt\eps.
\end{align}\end{subequations}
Therefore, for a subsequence, we have
\begin{subequations}\label{conv-2}\begin{align}\label{conv-u-2}
&u_\eps\to u&&\text{weakly* in}\ W^{1,\infty}(I;L^2(\Omega;\R^3))\cap 
H^1(I;H^1(\Omega;\R^3)),
\\&\pi_\eps\to\pi&&\text{weakly in}\ H^1(I;L^2(\Omega;\Rsym)),
\\&\label{conv-pi-L}
\pi_\eps|_{\QF}^{}\!\to0&&\text{strongly in }\ H^1(I;L^2(\OF;\Rsym)),
\\&\label{conv-a}
\alpha_\eps|_{\QS}^{}\!\to\alpha&&\text{weakly* in}\ L^\infty(I;H^1(\OS))\cap H^1(I;L^2(\OS)),
\\&\label{conv-phi}
\phi_\eps\to\phi&&\text{weakly* in}\  L^\infty(I;H^1(\R^3)),\ \text{ and also}
\\&\label{conv-el-2}
e_{\rm el,\eps}\to e_{\rm el}&&\text{strongly in }L^2(Q;\Rsym).
\end{align}\end{subequations} 
For \eq{conv-pi-L}, we used $\DT\pi_\eps\to0$ strongly in 
$L^2(\QF;\R^{3\times 3})$ due to \eq{est-pi-2} together with the assumption 
$\pi_0=0$ on $\OF$.
For \eq{conv-el-2}, like \eq{KV-damage-strong-e(t)} when taking into
account that $\pi=0$ on $\QF$, we have
\begin{align}\nonumber
&\!\int_\Omega\frac12\bbD_{_{\rm KV}}(e_{{\rm el},\eps}(t){-}e_{\rm el}(t))\Colon(e_{{\rm el},\eps}(t){-}e_{\rm el}(t))\,
\d x
\\[-.4em]\nonumber
&\!\le\int_\Omega\frac12\bbD_{_{\rm KV}}(e_{{\rm el},\eps}(t){-}e_{\rm el}(t))\Colon(e_{{\rm el},\eps}(t){-}e_{\rm el}(t))\,\d x
\\[-.4em]&\nonumber\qquad\qquad
+\int_{\OS}\!\!\frac12\bbD_{_{\rm MX}}(\pi_\eps(t){-}\pi(t))\Colon(\pi_\eps(t){-}\pi(t))
\,\d x+\int_{\OF}\!\!\frac1{2\eps}\bbD_{_{\rm MX}}\pi_\eps(t)\Colon\pi_\eps(t)\,\d x
\\[-.1em]&\nonumber
+\!\int_0^t\!\!\int_\Omega\!\!K_{_{\rm E}}|{\rm sph}(e_{{\rm el},\eps}{-}e_{\rm el})|^2\!
+\big(\pl_e^{}\calG_{_{\rm E},\eps}(\alpha_\eps,e_{{\rm el},\eps})
{-}\pl_e^{}\calG_{_{\rm E},\eps}(\alpha_\eps,e_{\rm el})\big)
\Colon(e_{{\rm el},\eps}{-}e_{\rm el})\,\d x\d t
\\&\nonumber=\int_0^t\!\!\bigg(\int_\Omega-(\varrho\DDT u_\eps+f_\eps)\Cdot(u_\eps{-}u)
-\big(\bbD_{_{\rm KV}}\DT e_{\rm el}+K_{_{\rm E}}{\rm sph}e_{\rm el}
+\pl_e^{}\calG_{_{\rm E},\eps}(\alpha_\eps,e_{\rm el})\big)\Colon
\\[-.4em]&\nonumber\qquad\qquad\qquad\qquad\qquad\qquad\qquad
\Colon(e_{{\rm el},\eps}{-}e_{\rm el})\,\d x
-\int_{\OS}\!\!\bbD_{_{\rm MX}}\DT\pi\Colon(\pi_\eps{-}\pi)\,\d x\!\bigg)\,\d t
\\&\nonumber=
\int_0^t\!\!\bigg(\int_\Omega\!\varrho\DT u_\eps\Cdot(\DT u_\eps{-}\DT u)
-f_\eps\Cdot(u_\eps{-}u)
+\Big(\pl_e^{}\calG_{_{\rm E},\eps}(\alpha,e_{\rm el}){-}\pl_e^{}\calG_{_{\rm E},\eps}(\alpha_\eps,e_{\rm el})
\\[-.4em]&\nonumber\qquad\qquad\qquad\ \ -
\bbD_{_{\rm KV}}\DT e_{\rm el}-K_{_{\rm E}}{\rm sph}e_{\rm el}
-\pl_e^{}\calG_{_{\rm E},\eps}(\alpha,e_{\rm el})\Big)\Colon(e_{{\rm el},\eps}{-}e_{\rm el})\,\d x
\\[-.2em]&\qquad\qquad
-\int_{\OS}\!\!\bbD_{_{\rm MX}}\DT\pi\Colon(\pi_\eps{-}\pi)\,\d x\!\bigg)\,\d t
-\int_\Omega\varrho\DT u_\eps(t)\Cdot(u_\eps(t){-}u(t))\,\d x\ \to\ 0\,
\label{KV-damage-strong-e(t)-2}
\end{align}
with $f_\eps=\varrho(\omega{\times}(\omega{\times}(x{+}u_\eps))+
2\omega{\times}\DT u_\eps+\nabla\phi_\eps)$
bounded in $L^2(Q;\R^3)$ and where we now used 
the continuity of the Nemytski\u\i\ mapping induced by 
$\pl_e^{}\calG_{_{\rm E},\eps}(\cdot,e_{\rm el})=\pl_e^{}\calG_{_{\rm E}}(\cdot,e_{\rm el})$ 
on $\QS$ while, on $\QF$, we use that simply 
$\pl_e^{}\calG_{_{\rm E},\eps}(\alpha_\eps,e_{\rm el})=
\eps\pl_e^{}\calG_{_{\rm E}}(\alpha_\eps,e_{\rm el})\to0$ strongly in 
$L^2(\QF)$; recall the scaling \eq{G-K-eps}.

The limit passage towards the weak solution in the sense of 
Definition~\ref{def-2} is then simple by the continuity with 
respect to the convergences \eq{conv-2}. More details deserve 
only the limit passage in 
\eq{system-alpha-w} which, written for $(u_\eps,\pi_\eps,\alpha_\eps)$ and
omitting the terms $\eps\int_{\QS}\!\eta(\DT\alpha_\eps)\,\d x\d t\ge0$ and
$\eps\int_{\OF}\frac\kappa2|\nabla\alpha_\eps(T)|^2\,\d x\ge0$,
reads as:
\begin{align}
&\nonumber\!\!\!\!\!\!\int_{\QS}\!\!\!\big(v{-}\DT\alpha_\eps\big)
\partial_\alpha\calG_{_{\rm E},\eps}(\alpha_\eps,{\rm dev}(e(u_\eps){-}\pi_\eps))
+\kappa\nabla\alpha_\eps\Cdot\nabla v+\eta(v)\,\d x\d t
+\int_{\OS}\!\!\frac\kappa2|\nabla\alpha_0|^2\,\d x
\\&\qquad\nonumber
+\eps\bigg(\int_{\QF}\!\!\!\big(v{-}\DT\alpha_\eps\big)
\partial_\alpha\calG_{_{\rm E},\eps}(\alpha_\eps,{\rm dev}(e(u_\eps){-}\pi_\eps))
+\kappa\nabla\alpha_\eps\Cdot\nabla v+\eta(v)\,\d x\d t
\\&\qquad\qquad\quad
+\int_{\OF}\!\!\frac\kappa2|\nabla\alpha_0|^2\,\d x\bigg)
\ge\int_{\QS}\!\eta(\DT\alpha_\eps)\,\d x\d t
+\int_{\OS}\frac\kappa2|\nabla\alpha_\eps(T)|^2\,\d x.
\label{system-alpha-w-2}
\end{align}
Now we use 
$$
\bigg|\eps\!\int_{\QF}\!\!\!\!\DT\alpha_\eps\partial_\alpha\calG_{_{\rm E},\eps}(\alpha_\eps,
e_{\rm el,\eps})
\d x\d t\bigg|
\le\eps C\|\DT\alpha_\eps\|_{L^2(\QF)}^{}\|1{+}|e_{\rm el,\eps}|\,\|_{L^2(Q)}^{}
=\mathscr{O}(\sqrt\eps)\to0
$$ 
due to \eq{est-DT-a-2}, and
$$
\bigg|\eps\!\int_{\QF}\!\!\!\!\kappa\nabla\alpha_\eps\Cdot\nabla v\d x\d t\bigg|
\le\eps\big(\sup_\Omega\kappa\big))\|\nabla\alpha_\eps\|_{L^2(\QF;\R^3)}^{}\|\nabla v\|_{L^2(Q;\R^3)}^{}
=\mathscr{O}(\sqrt\eps)\to0
$$ 
due to \eq{est-nabla-a-2}.
In the limit we thus obtain \eq{system-alpha-w} on $\OS$.

Altogether, we proved that the limit $(u,\pi,\alpha,\phi)$
solves the initial-boundary-value problem
\eq{system++}--{\rm(\ref{BC-IC}a,b)}--\eq{BC-IC-alpha+}
in the sense of Definition~\ref{def-2}.

In addition, this solution conserves energy. This can be shown
again by using that $\sqrt\varrho\DDT u$ in duality with 
$\sqrt\varrho\DT u$ so that \eq{calculus1} holds, and 
that ${\rm div}(\kappa\nabla\alpha)\in L^2(\QS)$ and 
$\DT\alpha\in L^2(\QS)$ so that also \eq{calculus2} holds
but now with $\QS$ and $\OS$ instead of $Q$ and $\Omega$, respectively.

To prove \eq{Max-limit}, we use \eq{engr-int} written for 
$(u_\eps,\pi_\eps,\alpha_\eps,\phi_\eps)$ with $t=T$ and the scaling 
\eq{BC-IC-alpha+}, namely
\begin{align}\nonumber
&\limsup_{\eps\to0}\int_{\QF}\frac1\eps\bbD_{_{\rm MX}}\DT\pi_\eps\Colon\DT\pi_\eps
+\eps\DT\alpha_\eps\zeta(\DT\alpha_\eps)\,\d x\d t
\\\nonumber&\hspace{0em}
\le\lim_{\eps\to0}\int_{\R^3}\!\!\varrho_{\rm ext}(T)\phi_\eps(T)\,\d x
-\int_\Omega\!\varrho\nabla\phi_\eps(T)\Cdot u_\eps(T)-\varrho\phi_\eps(T)
\,\d x
-\!\int_0^T\!\!\!\int_{\R^3}\!\DT\varrho_{\rm ext}\phi_\eps\,\d x\d t
\\[-.4em]&\nonumber\hspace{6em}
-\int_Q(\omega{\times}(\omega{\times}(x{+}u_\eps)\cdot\DT u_\eps\,\d x\d t
-\int_{\OF}\!\!\eps\calG_{_{\rm E}}(\alpha_\eps(T),e_{\rm el,\eps}(T))\,\d x
\\[-.4em]&\nonumber\hspace{0em}
-\liminf_{\eps\to0}\bigg(\int_\Omega\frac\varrho2|\DT u_\eps(T)|^2
+\frac32 K_{_{\rm E}}|{\rm sph}\,e_{\rm el,\eps}(T)|^2
\,\d x
\\[-.3em]\nonumber
&\hspace{5em}+\int_{\OS}\!\!\calG_{_{\rm E}}(\alpha_\eps(T),e_{\rm el,\eps}(T))
+\frac\kappa2|\nabla\alpha_\eps(T)|^2\,\d x
+
\int_{\R^3}\frac1{
8\uppi g}|\nabla\phi_\eps(T)|^2\,\d x
\\[-.3em]&\hspace{6em}\nonumber
+\int_Q\!\bbD_{_{\rm KV}}\DT e_{\rm el,\eps}\Colon \DT e_{\rm el,\eps}
\,\d x\d t+\int_{\QS}\bbD_{_{\rm MX}}\DT\pi_\eps\Colon\DT\pi_\eps
+\DT\alpha_\eps\eta'(\DT\alpha_\eps)\,\d x\d t\bigg)
+E_0
\\[-.4em]&\nonumber\hspace{0em}
\le\int_{\R^3}\!\!\varrho_{\rm ext}(T)\phi(T)-\frac1{
8\uppi g}|\nabla\phi(T)|^2\,\d x
-\int_\Omega\!\varrho\nabla\phi(T)\Cdot u(T)-\varrho\phi(T)
\,\d x
\\[-.4em]&\nonumber\hspace{1em}
-\!\int_0^T\!\!\!\int_{\R^3}\!\DT\varrho_{\rm ext}\phi\,\d x\d t
-\int_Q\bbD_{_{\rm KV}}\DT e_{\rm el}\Colon \DT e_{\rm el}+(\omega{\times}(\omega{\times}(x{+}u)\cdot\DT u\,\d x\d t
\\[-.1em]&\nonumber\hspace{1em}
-\int_\Omega\frac\varrho2|\DT u(T)|^2
+\frac32 K_{_{\rm E}}|{\rm sph}\,e_{\rm el}(T)|^2
\,\d x
-\!\int_{\OS}\!\!\!\calG_{_{\rm E}}(\alpha(T),e_{\rm el}(T))+\frac\kappa2|\nabla\alpha(T)|^2\,\d x
\\[-.3em]
&\hspace{13em}
-\!\int_{\QS}\!\!\!\bbD_{_{\rm MX}}\DT\pi\Colon\DT\pi
+\DT\alpha\eta'(\DT\alpha)\,\d x\d t+E_0=0,
\label{engr-int-eps}\end{align}
where we again used the notation \eq{C-D} and $E_0$ from 
\eq{def-E0} now with $\kappa=0$ on $\OF$. The 
first inequality arises from ``forgetting'' the 
term $\frac\kappa2|\nabla\alpha_\eps(T)|^2$ on $\OF$, while the
second inequality is by weak lower semicontinuity.
The last equality in \eq{engr-int-eps} expresses the energy conservation 
for the limit system, discussed already above. 
%
Thus \eq{Max-limit} is proved.
$\Box$

\begin{remark}[{\sl Strong convergence of total strains}]\upshape
From \eq{KV-damage-strong-e(t)-2}, one can even see that also $\pi_\eps\to\pi$
strongly so that, together with \eq{conv-el-2}, even the total strain
$e(u_\eps)=e_{{\rm el},\eps}+\pi_\eps$ converges strongly in $L^2(Q;\Rsym)$. We 
did not need this additional result in the above convergence proof, however.
\end{remark}

\begin{remark}[{\sl Damage flow rule almost everywhere}]\label{rem-flow-rule}
\upshape
 Since $\partial\zeta(\DT\alpha)$ is bounded in $L^2(Q)$ in our model,
by comparison, we have also 
${\rm div}(\kappa\nabla\alpha)\in\pl\eta(\DT\alpha)
+\partial_\alpha\calG_{_{\rm E}}(\alpha,{\rm dev}\,e_{\rm el})$ bounded in $L^2(Q)$, cf.\
\eq{system-alpha++}. Since also $\DT\alpha\In L^2(Q)$, 
the formula $\int_Q\DT\alpha{\rm div}(\kappa\nabla\alpha)\,\d x\d t
=\frac12\int_\Omega\kappa|\nabla\alpha_0|^2-\kappa|\nabla\alpha(T)|^2\,\d x$
rigorously holds, and we can write \eq{system-alpha-w} in
as the original inequality \eq{system-alpha-w-} holding a.e.\ on $Q$.
The integral form \eq{system-alpha-w} is however suitable for the limit
passages, in contrast to \eq{system-alpha-w-}.
\end{remark}


\begin{definition}\label{def-3}
Assuming, beside \eq{ass}, also 
$e_{\rm el,0}=e(u_0){-}\pi_0\in H^1(\Omega;\R^{3\times3})$ and 
$v_0\in H^2(\Omega;\R^3)$, then the weak solution to the system 
{\rm(\ref{system++}a,d-f)}--{\rm(\ref{BC-IC}a,b)}--\eq{BC-IC-alpha+}
with \eq{system+++} is understood as a five-tuple
$(u,\sigma_{\rm sph},\sigma_{\rm dev},\pi,\phi)\in W^{1,\infty}(I;L^2(\Omega;\R^3))\times L^2(Q;\Rsph)\times L^2(Q;\Rdev)\times H^1(I;L^2(\OS;\Rsym))\times W^{1,\infty}(I;H^1(\R^3))$
if the integral identity
\begin{align}
\int_{\QS}\!\!\!\sigma_{\rm dev}\Colon{\rm dev}\,e(v)\,\d x\d t+
\int_Q\!
\sigma_{\rm sph}\Colon{\rm sph}\,e(v)-\r\DT{u}\Cdot\DT{v}+f\Cdot v\,\d x\d t
=\!\int_\Omega\!\r v_0\Cdot v(0,\cdot)\,\d x
\end{align}
with $f$ from \eq{system-u+} holds for any $v\In H^1(Q;\R^3)$ with $v(T)=0$, 
further \eq{system-pi++} hold a.e.\ on $\QS$ and \eq{system+++} relating 
$(\sigma_{\rm dev},\sigma_{\rm sph})$ with $(u,\pi)$ hold a.e.\ on $Q$, and also 
the initial conditions $u(0,\cdot)=u_0$ and $\pi(0,\cdot)=\pi_0$ hold 
a.e.\ on $\Omega$.
\end{definition}

Note that, controlling $\sigma_{\rm sph}$ and $\sigma_{\rm dev}$ in the 
above definition, we have implicitly also included the information 
$$
u|_{\QS}\!\In H^1(I;H^1(\Omega;\R^3)) \quad\text{ and }\quad  
u\In L^2(I;L_{\rm div}^2(\Omega;\R^3))
$$
with $L_{\rm div}^2(\Omega;\R^3)
=\{u\In L^2(\Omega;\R^3);\ {\rm div}\,u\In L^2(\Omega)\}$,
while the strain $e(u)$ is not defined in the fluidic regions $\QF$.
Also one has $e(u){-}\pi\In H^1(I;L^2(\OS;\Rsym))$.

\medskip

\noindent{\it Sketched proof of Proposition~\ref{prop-3}.}
We first prove a certain regularity by differentiating in time the system
(\ref{system++}a,d) with \eq{system+++} written for the solution obtained in 
Proposition~\ref{prop-1} and denoted now as $(u_\eps,\pi_\eps,\alpha_\eps,\phi_\eps)$, 
and employing the test by $(\DDT u_\eps,\DDT\pi_\eps)$. Taking into account the 
scaling \eq{KV-eps} and using again the orthogonality for the Coriolis force 
(now in terms of accelerations) as $(2\varrho\omega{\times}\DDT u)\cdot\DDT u=0$, 
this gives
\begin{align}\nonumber
&\int_\Omega\frac\varrho2|\DDT u_\eps(t)|^2
+\frac32K_{_{\rm E}}|{\rm sph}\,\DT e_{\rm el,\eps}(t)|^2\,\d x
\\[-.4em]&\nonumber\qquad
+\int_0^t\!\!\int_{\OS}\!\!\!\bbD_{_{\rm KV}}\DDT e_{\rm el,\eps}\Colon\DDT e_{\rm el,\eps}
+\bbD_{_{\rm MX}}\DDT\pi_\eps\Colon\DDT \pi_\eps\,\d x\d t
+\eps\!\int_0^t\!\!\int_{\OF}\!\!\!\bbD_{_{\rm KV}}e(\DDT u_\eps)\Colon e(\DDT u_\eps)\,\d x\d t
\\&\nonumber
=\!\int_\Omega\frac\varrho2|\DDT u_\eps(0)|^2
+\frac32K_{_{\rm E}}|{\rm sph}\,\DT e_{\rm el,\eps}(0)|^2\d x
-\int_0^t\!\!\int_\Omega\!
\varrho\big(\omega{\times}(\omega{\times}(x{+}\DT u_\eps)){+}
\nabla\DT\phi_\eps\big)\Cdot \DDT u_\eps
\\[-.4em]&\nonumber\qquad
-\int_0^t\!\!\int_{\OS}\!\!\!\Big(\pl_{ee}^2\calG_{_{\rm E}}(\alpha_\eps,e_{\rm el,\eps})\Colon
{\rm dev}\,\DT e_{\rm el,\eps}\!+
\pl_{e\alpha}^2\calG_{_{\rm E}}(\alpha_\eps,e_{\rm el,\eps})\DT\alpha_\eps\Big)
\Colon {\rm dev}\,\DDT e_{\rm el,\eps}\d x\d t
\\&\nonumber
\le\frac12\|\varrho\|_{L^\infty(\Omega)}^{}\|\DDT u_\eps(0)\|_{L^2(\Omega;\R^3)}^2
+\frac32K_{_{\rm E}}\|{\rm sph}(e(v_0){-}\DT\pi_0)\|_{L^2(\Omega;\R^{3\times 3})}^2
\\&\nonumber\quad
+\|\sqrt\varrho\|_{L^\infty(\Omega)}^{}
\big\|\omega{\times}(\omega{\times}(x{+}\DT u_\eps)){+}
\nabla\DT\phi_\eps\big\|_{L^2(Q;\R^3)}^2+\int_0^t\!\!\int_\Omega\frac\varrho2|\DDT u_\eps(t)|^2\,\d x\d t
\\&\nonumber\quad
+\frac1c\Big(
\|\pl_{ee}^2\calG_{_{\rm E}}(\alpha_\eps,e_{\rm el,\eps})\|_{L^\infty(\QS;\R^{3^4})}^2
\|e(\DT u_\eps)\|_{L^2(Q;\R^{3\times 3})}^2
\\&
\quad+
\|\pl_{e\alpha}^2\calG_{_{\rm E}}(\alpha_\eps,e_{\rm el,\eps})\|_{L^\infty(Q;\R^{3\times 3})}^2
\|\DT\alpha_\eps\|_{L^2(Q)}^2\Big)
+\epsilon
\|\DDT e_{\rm el,\eps}\|_{L^2([0,t]\times\OS;\R^{3\times 3})}^2
\label{est-regular}
\end{align}
with $c:=\frac12\|\bbD_{_{\rm KV}}\|_{L^\infty(\Omega;\R^{3^4})}$; 
then the last term in \eq{est-regular} can be absorbed in the left
hand side. Note that we needed the assumptions \eq{ass-growth-G} to have
$\pl_{e\alpha}^2\calG_{_{\rm E}}(\alpha_\eps,e_{\rm el,\eps})$ and
$\pl_{ee}^2\calG_{_{\rm E}}(\alpha_\eps,e_{\rm el,\eps})$ apriori bounded.
We further use
\begin{align}
\DDT u_\eps(0)=\frac1\varrho\Big({\rm div}\,\sigma_{0,\eps}
+\varrho(\omega{\times}(\omega{\times}(x{+}u_0))+
2\omega{\times}v_0+\nabla\phi_0)
\Big)\in L^2(\Omega;\R^3)
\label{est-regular+}
\end{align}
with the initial stress $\sigma_{0,\eps}=\sigma_{\rm sph}+\sigma_{\rm dev}$ 
with $\sigma_{\rm sph}$ and $\sigma_{\rm dev}$ from (\ref{system++}b,c)
with $(K_{_{\rm KV}},G_{_{\rm KV}})=(K_{_{\rm KV},\eps},G_{_{\rm KV},\eps})$ from \eq{KV-eps},
and with $e_{\rm el}=e(u_0)-\pi_0$ and $\DT e_{\rm el}=e(v_0)-\DT\pi_0$,
and $\phi_0\in H^1(\R^3)$ solving 
${\rm div}(\nabla\phi_0/(4\uppi g)+\varrho u_0)=\varrho+\varrho_{\rm ext}(0)$.
Note that indeed $\DDT u_\eps(0)\in L^2(\Omega;\R^3)$
due to the assumptions 
$e_{\rm el,0}=e(u_0){-}\pi_0\in H^1(\Omega;\R^{3\times3})$ and 
$v_0\in H^2(\Omega;\R^3)$ involved in Definition~\ref{def-3}. 
Also note that $\DT\pi_0:=\DT\pi(0)$ is involved in \eq{est-regular}
and implicitly also in \eq{est-regular+}, although we do not 
prescribe any initial condition on $\DT\pi$. Yet, from \eq{system-pi++},
we can read $\DT\pi(0)=\bbD_{_{\rm MX}}^{-1}\sigma_{0}$ on $\OS$ while 
$\DT\pi(0)=0$ on $\OF$ because we have now $\bbD_{_{\rm MX}}^{-1}=0$ on $\OF$. 
To estimate $\sigma_{0}$ let us realize that
$\|\sigma_{0}\|_{H^1(\Omega;\R^{3\times3})}\le
\|\bbD_{_{\rm KV}}(e(v_0){-}\DT\pi(0))
+\pl_e\calG_{{\rm E}}(\alpha_0,e_{\rm el,0})\|_{H^1(\Omega;\R^{3\times3})}\le
C+\|\bbD_{_{\rm KV}}\bbD_{_{\rm MX}}^{-1}\sigma_{0}\|_{H^1(\Omega;\R^{3\times3})}$
with some $C$ depending on $\|v_0\|_{H^2(\Omega;\R^3)}$, $\|\alpha_0\|_{H^1(\Omega)}$,
and $\|e_{\rm el,0}\|_{H^1(\Omega;\R^{3\times3})}$; note that these quantities 
have been supposed finite in Definition~\ref{def-3}.
Thus, using \eq{ass-small-KV/MX}, we get the desired bound on $\sigma_{0}$.
Hence, $\DT e_{\rm el}(0)$ which occurs in \eq{est-regular}
itself can be estimated.

Using Gronwall's inequality for \eq{est-regular} then yields the estimates 
\begin{subequations}\label{est-3}\begin{align}\label{est-u-3}
&\|u_\eps\|_{W^{2,\infty}(I;L^2(\Omega;\R^3))\,\cap\,H^2(I;H^1(\OS;\R^3))}^{}\le C,
\\&\label{est-pi-3}
\|\DDT\pi_\eps\|_{L^2(\QS;\R^{3\times 3})}^{}\le C,\ \ \text{ and }\ \ 
\\&\label{est-pi-3+}
\|u_\eps\|_{H^2(I;H^1(\OF;\R^3))}^{}\le C/\sqrt\eps,
\intertext{which are now at disposal together with \eq{est-phi} and 
\eq{est-a-2}.
Also we have}
&\|\sigma_\eps\|_{L^2(Q;\Rsym)}^{}\le C\ \ \text{ with }\ \sigma_\eps=\begin{cases}
\eps\bbD_{_{\rm KV}}e(\DT u_\eps)+3K_{_{\rm KV}}{\rm sph}\,e(u_\eps)\!\!
&\text{on }\QF,\\\ 
\ \bbD_{_{\rm KV}}\DT e_{\rm el,\eps}+\bbC_{_{\rm E}}(\alpha_\eps,e_{\rm el,\eps})
&\text{on }\QS,
\end{cases}\end{align}\end{subequations}
where we again use the notation \eq{C-D} and $e_{\rm el,\eps}=e(u_\eps)-\pi_\eps$.

By Banach's selection principle, we consider a weakly* convergent subsequence
respecting the topologies specified in (\ref{est-3}a,b) and \eq{est-a-2}
together with 
the convergence \eq{conv-phi} and also $\sigma_\eps\to\sigma$ weakly in 
$L^2(Q;\Rsym)$. Then we put naturally $\sigma_{\rm sph}:={\rm sph}\,\sigma$ and 
$\sigma_{\rm dev}:={\rm dev}\,\sigma$. 
As $\|e(\DT u_\eps)\|_{L^2(\QF;\R^{3\times3})}=\mathscr{O}(1/\sqrt\eps)$,
we have $\|\eps\bbD_{_{\rm KV}}e(\DT u_\eps)\|_{L^2(\QF;\R^{3\times3})}
=\mathscr{O}(\sqrt\eps)\to0$ so that
$\sigma=3K_{_{\rm KV}}{\rm sph}\,e(u)$ on $\QF$.
To identify $\sigma=\bbD_{_{\rm KV}}\DT e_{\rm el}+\bbC_{_{\rm E}}(\alpha_\eps,e_{\rm el})$
in the solid regions $\QS$,
due to nonlinearities $\bbC_{_{\rm E}}(\alpha,\cdot)$, we again need to prove 
$e_{\rm el,\eps}\to e_{\rm el}$ strongly in 
$L^2(\QS;\Rsym)$. To this goal, one is to modify 
\eq{KV-damage-strong-e(t)-2} to be used on $\OS$ rather than $\Omega$.
The peculiarity is that $u_\eps-u$ is no longer a legitimate
test function because $e(u)$ is not well defined in the fluidic regions 
$\QF$. For this reason, we take some smooth $\tilde u_\eps$ that will
approximate $u$ strongly in $L^\infty(I;L_{\rm div}^2(\Omega;\R^3))$ 
and even $\tilde u_\eps|_{\QS}\!\to u|_{\QS}$ strongly in 
$H^1(I;H^1(\Omega;\R^3))$, and we can assume that this convergence is 
sufficiently slow so that 
\begin{align}\label{approximation}
\|\tilde u_\eps\|_{L^2(I;H^1(\OF;\R^3))}\le 1/\sqrt[4]{\eps}.
\end{align}
Then, denoting $\tilde e_{{\rm el},\eps}=e(\tilde u_\eps)-\pi$,
we have $\tilde e_{{\rm el},\eps}(t)\to e_{\rm el}(t)$ for a.a.\ $t\In I$
and we can write
\begin{align}\nonumber
&\limsup_{\eps\to0}\!\int_{\OS}\frac14\bbD_{_{\rm KV}}(e_{{\rm el},\eps}(t){-}e_{\rm el}(t))\Colon(e_{{\rm el},\eps}(t){-}e_{\rm el}(t))\,
\d x
\\[-.4em]\nonumber
&\le\limsup_{\eps\to0}\!\int_{\OS}\frac12\bbD_{_{\rm KV}}(e_{{\rm el},\eps}(t){-}\tilde e_{{\rm el},\eps}(t))\Colon(\tilde e_{{\rm el},\eps}(t){-}e_{{\rm el},\eps}(t))\,\d x
\\[-.4em]&\nonumber\qquad\qquad\qquad
+\lim_{\eps\to0}\!\int_{\OS}\frac12\bbD_{_{\rm KV}}(\tilde e_{{\rm el},\eps}(t){-}e_{\rm el}(t))\Colon(\tilde e_{{\rm el},\eps}(t){-}e_{\rm el}(t))\,\d x
\displaybreak\\[-.4em]\nonumber
&\!\le\limsup_{\eps\to0}\bigg(\!\int_{\OS}\!\!\bigg(\frac12\bbD_{_{\rm KV}}(e_{{\rm el},\eps}(t){-}\tilde e_{{\rm el},\eps}(t))\Colon(e_{{\rm el},\eps}(t){-}\tilde e_{{\rm el},\eps}(t))
\\[-.4em]&\nonumber\qquad
+
\frac12\bbD_{_{\rm MX}}(\pi_\eps(t){-}\pi(t))\Colon(\pi_\eps(t){-}\pi(t))\bigg)
\,\d x
+\!\int_0^t\!\!\int_\Omega\!\!K_{_{\rm E}}|{\rm sph}(e_{{\rm el},\eps}{-}\tilde e_{{\rm el},\eps})|^2\,\d x\d t
\\[-.1em]&\nonumber\qquad\qquad
+\!\int_0^t\!\!\int_{\OS}\!\!\!\big(\pl_e^{}\calG_{_{\rm E}}(\alpha_\eps,e_{{\rm el},\eps})
{-}\pl_e^{}\calG_{_{\rm E}}(\alpha_\eps,\tilde e_{{\rm el},\eps})\big)
\Colon(e_{{\rm el},\eps}{-}\tilde e_{{\rm el},\eps})\,\d x\d t\bigg)
\\&\nonumber\le\lim_{\eps\to0}\bigg(
\int_{\OF}\!\!\!\frac\eps2\bbD_{_{\rm KV}}e(u_0)\,\d x
-\!\int_0^t\!\!\int_{\OF}\!\!\!\eps\bbD_{_{\rm KV}}e(\DT u_\eps)\Colon e(\tilde u_\eps)
\,\d x\d t
\\&\nonumber\qquad
-
\int_0^t\!\!\int_\Omega\!(\varrho\DDT u_\eps+f_\eps)\Cdot(u_\eps{-}\tilde u_\eps)
-\!\int_0^t\!\!\int_{\OS}\!\!\!\bbD_{_{\rm MX}}\DT\pi\Colon(\pi_\eps{-}\pi)
\\[-.4em]&
\qquad
+\big(\bbD_{_{\rm KV}}\DT{\tilde e}_{{\rm el},\eps}\!
+K_{_{\rm E}}{\rm sph}\tilde e_{{\rm el},\eps}\!
+\pl_e^{}\calG_{_{\rm E}}(\alpha_\eps,\tilde e_{{\rm el},\eps})\big)
\Colon(e_{{\rm el},\eps}{-}\tilde e_{{\rm el},\eps})
\,\d x\d t\bigg)=0,
\label{KV-damage-strong-e(t)-3}
\end{align}
where the third inequality have arisen by ``forgetting'' the nonnegative
term $\int_{\OF}\!\!\!\eps\bbD_{_{\rm KV}}e(u_\eps(T))\Colon e(u_\eps(T))\,\d x$.
Note the we used \eq{approximation} together with \eq{est-pi-3+} for the
estimate
\begin{align}\nonumber
&\Big|\!\int_0^t\!\!\int_{\OF}\!\!\!\eps\bbD_{_{\rm KV}}e(\DT u_\eps)\Colon e(\tilde u_\eps)
\,\d x\d t\Big|\le\eps\|\bbD_{_{\rm KV}}\|_{L^\infty(\Omega;\R^{3^4})}\times
\\[-.6em]&\qquad\qquad
\times\|e(\DT u_\eps)\|_{L^2(\QF;\R^{3\times3})}\|e(\tilde u_\eps)\|_{L^2(\QF;\R^{3\times3}))}
=\eps\mathscr{O}\Big(\frac1{\sqrt[3/4]{\eps}}\Big)
=\mathscr{O}(\sqrt[4]{\eps})\to0.
\nonumber
\end{align}
%
From \eq{KV-damage-strong-e(t)-3}, we thus have 
$e_{\rm el,\eps}|_{\QS}(t)\to e_{\rm el}|_{\QS}(t)$ at a.a. $t\In I$. Then, 
instead of \eq{conv-el-2}, by the Lebesgue theorem, we now proved
\begin{align}
\label{conv-el-2+}
e_{\rm el,\eps}|_{\QS}\to e_{\rm el}|_{\QS}&&\text{strongly in }L^2(\QS;\Rsym),
\end{align}
which is to be used for the limit passage in the nonlinear term
$\pl_e^{}\calG_{_{\rm E}}$.

The energy conservation now holds due to the proved regularity, as 
$\sqrt\varrho\DDT u\in L^\infty(I;L^2(\Omega;\R^3))$ is surely in duality with 
$\sqrt\varrho\DT u\in W^{1,\infty}(I;L^2(\Omega;\R^3))$, cf.\ \eq{est-u-3}.
Eventually, like we did in \eq{engr-int-eps}, we now can show that
\begin{align}\nonumber
&\limsup_{\eps\to0}\int_{\QF}\!\!\eps\bbD_{_{\rm KV}} e(\DT u_\eps)\Colon e(\DT u_\eps)\,\d x\d t\le\int_{\R^3}\!\!\varrho_{\rm ext}(T)\phi(T)-\frac1{
8\uppi g}|\nabla\phi(T)|^2\,\d x
\\[-.4em]&\nonumber\hspace{.5em}
-\int_\Omega\!\varrho\nabla\phi(T)\Cdot u(T)-\varrho\phi(T)
\,\d x-\!\int_0^T\!\!\!\int_{\R^3}\!\DT\varrho_{\rm ext}\phi\,\d x\d t
-\int_Q(\omega{\times}(\omega{\times}(x{+}u)\cdot\DT u\,\d x\d t
\\[-.4em]&\nonumber\hspace{1.5em}
-\int_\Omega\frac\varrho2|\DT u(T)|^2
+\frac32 K_{_{\rm E}}|{\rm sph}\,e_{\rm el}(T)|^2
\,\d x-\!\int_{\OS}\!\!\!\calG_{_{\rm E}}(\alpha(T),e_{\rm el}(T))+\frac\kappa2|\nabla\alpha(T)|^2\,\d x
\\[-.3em]
&\hspace{2.5em}
-\!\int_{\QS}\!\!\!\bbD_{_{\rm KV}}\DT e_{\rm el}\Colon \DT e_{\rm el}+\bbD_{_{\rm MX}}\DT\pi\Colon\DT\pi
+\DT\alpha\eta'(\DT\alpha)\,\d x\d t+E_0=0,
\label{engr-int-eps-2}
\end{align}
with $E_0$ from \eq{def-E0}. The last equality is due to the mentioned energy 
conservation. Thus \eq{KV-limit} is proved.
$\Box$

\begin{remark}[{\sl The successive convergence}]\upshape
Note that, in \eq{est-regular}, we benefited from having 
$G_{_{\rm E}}$ already pushed to zero on $\OF$ because the
viscosity on $\OF$ is (intentionally) not uniformly controlled,
being limited to zero. Thus the direct merging of both 
limit processes in Propositions~\ref{prop-2} 
and \ref{prop-3} is not possible. Of course, a suitable
scaling between these two would facilitate such connection.
\end{remark}

\begin{remark}[{\sl Discontinuities in $\varrho$ accros the Earth interfaces}]\upshape
The interfaces between the ocean beds and the mantle as well as the 
Gutenberg's discontinuity and  between the inner and outer core regions 
typically also exhibit discontinuities in
mass density $\varrho$, which is incompatible with the assumption 
$\varrho\in W^{1,3}(\Omega)$ in \eq{ass-rho} used in \eq{est-of-DDTu}.
Note that the additional regularity of the initial conditions involved
in Definition~\ref{def-3} allowed us to avoid this restriction and
consider a general, possibly discontinuous $\varrho\in L^\infty(\Omega)$.
A respective modification of the proofs of 
Propositions~\ref{prop-1}--\ref{prop-2} would be possible, too.
\end{remark}

\bigskip\noindent{\it Acknowledgments.} 
The author is thankful to Katharina Brazda and Ctirad Matyska for
many fruitful discussions about the models. 
Moreover, deep thanks are also for hospitality and support of the 
Erwin Schr\"odinger Institute, Univ.\ Vienna, 
and for the partial support of the 
Czech Science Foundation projects 16-03823S and 17-04301S and also 
by the Austrian-Czech projects 16-34894L (FWF/CSF),
as well as through the institutional support RVO:\,61388998 (\v CR).

\end{sloppypar}

\begin{thebibliography}{99}

\bibitem{BeZi01DRRM}
Y.~Ben-Zion.
\newblock Dynamic ruptures in recent models of earthquake faults.
\newblock {\em J. Mech. Phys. Solids}, 49:2209--2244, 2001.

\bibitem{BZamp09SRRS}
Y.~Ben-Zion and J.-P. Ampuero.
\newblock Seismic radiation from regions sustaining material damage.
\newblock {\em Geophys. J. Int.}, 178:1351--1356, 2009.

\bibitem{Bog77HECV}
D.~Boger.
\newblock A highly elastic constant-viscosity fluid.
\newblock {\em J. Non-Newtonian Fluid Mechanics}, 3:87--91, 1977.

\bibitem{Braz17EGEG}
K.~Brazda.
\newblock {\em The elastic-gravitational equations in global seismology with
  low regularity}.
\newblock PhD thesis, Univ. Wien, 2017.

\bibitem{BrHoHo17VFEE}
K.~Brazda, M.~V. de~Hoop, and G.~Hoermann.
\newblock Variational formulation of the earth's elastic-gravitational
  deformations under low regularity conditions.
\newblock Preprint arXiv:1702.04741, 2017.

\bibitem{DahTro98TGS}
F.~A. Dahlen and J.~Tromp.
\newblock {\em Theoretical global seismology}.
\newblock Princetown Univ. Press, Princetown, NJ, 1998.

\bibitem{GreNag65GTEP}
A.~Green and P.~Naghdi.
\newblock A general theory of an elastic-plastic continuum.
\newblock {\em Arch. Rational Mech. Anal.}, 18:251--281, 1965.

\bibitem{HBAD09DERC}
R.~A. {Harris et al.}
\newblock The {SCEC/USGS} dynamic earthquake rupture code verification
  exercise.
\newblock {\em Seismological Res. Lett.}, 80:119--126, 2009.

\bibitem{HuAmHel14ERMW}
Y.~Huang, J.-P. Ampuero, and D.~V. Helmberger.
\newblock Earthquake ruptures modulated by waves in damaged fault zones.
\newblock {\em J. of Geophysical Research: Solid Earth}, B9:3133--3154, 2014.

\bibitem{KaLaAm08SEMS}
Y.~Kaneko, N.~Lapusta, and J.-P. Ampuero.
\newblock Spectral element modeling of spontaneous earthquake rupture on rate
  and state faults: Effect of velocity-strengthening friction at shallow
  depths.
\newblock {\em J. Geophysical Res.}, 113:B09317, 2008.

\bibitem{KomTro02SESG}
D.~Komatitsch and J.~Tromp.
\newblock Spectral-element simulations of global seismic wave propagation -
  {I}. validation.
\newblock {\em Geophys. J. Int.}, 149:390--412, 2002.

\bibitem{KomTro02SESG2}
D.~Komatitsch and J.~Tromp.
\newblock Spectral-element simulations of global seismic wave propagation -
  {II}. three-dimensional models, oceans, rotation and self-gravitation.
\newblock {\em Geophys. J. Int.}, 150:303--318, 2002.

\bibitem{KooDum11VEIC}
L.~Koot and M.~Dumberry.
\newblock Viscosity of the {E}arth's inner core: Constraints from nutation
  observations.
\newblock {\em Earth and Planetary Science Letters}, 308(3):343--349, 2011.

\bibitem{LayWal95MGS}
T.~Lay and T.~C. Wallace.
\newblock {\em Modern global seismology}.
\newblock Acad. Press, San Diego, 1995.

\bibitem{LyaBZ14DBRM}
V.~Lyakhovsky and Y.~Ben-Zion.
\newblock Damage-breakage rheology model and solid-granular transition near
  brittle instability.
\newblock {\em J. Mech. Phys. Solids}, 64:184--197, 2014.

\bibitem{LHAB09NDRW}
V.~Lyakhovsky, Y.~Hamiel, J.-P. Ampuero, and Y.~Ben-Zion.
\newblock Non-linear damage rheology and wave resonance in rocks.
\newblock {\em Geophys. J. Int.}, 178:910--920, 2009.

\bibitem{LyHaBZ11NLVE}
V.~Lyakhovsky, Y.~Hamiel, and Y.~Ben-Zion.
\newblock A non-local visco-elastic damage model and dynamic fracturing.
\newblock {\em J. Mech. Phys. Solids}, 59:1752--1776, 2011.

\bibitem{LyaMya84BECS}
V.~Lyakhovsky and V.~P. Myasnikov.
\newblock On the behavior of elastic cracked solid.
\newblock {\em Phys. Solid Earth}, 10:71--75, 1984.

\bibitem{MeaFur13SSWO}
T.~Maedae and T.~Furumura.
\newblock {FDM} simulation of seismic waves, ocean acoustic waves, and tsunamis
  based on tsunami-coupled equations of motion.
\newblock {\em Pure Appl. Geophys.}, 170:109--127, 2013.

\bibitem{PPAB12TDDR}
C.~Pelties, J.~{de la Puente}, J.-P. Ampuero, G.~B. Brietzke, and M.~K\"{a}ser.
\newblock Three-dimensional dynamic rupture simulation with a high-order
  discontinuous galerkin method on unstructured tetrahedral meshes.
\newblock {\em J. Geophys. Res.}, 117:B02309, 2012.

\bibitem{RajRou03EDSM}
K.~R. Rajagopal and T.~Roub{\'\i}{\v{c}}ek.
\newblock On the effect of dissipation in shape-memory alloys.
\newblock {\em Nonlinear Anal., Real World Appl.}, 4:581--597, 2003.

\bibitem{Roub13NPDE}
T.~Roub{\'\i}{\v{c}}ek.
\newblock {\em Nonlinear Partial Differential Equations with Applications}.
\newblock Birkh\"auser, Basel, 2nd edition, 2013.

\bibitem{Roub14NRSF}
T.~Roub{\'{\i}}{\v{c}}ek.
\newblock A note about the rate-and-state-dependent friction model in a
  thermodynamical framework of the {B}iot-type equation.
\newblock {\em Geophysical J. Intl.}, 199:286--295, 2014.

\bibitem{Roub17GMHF}
T.~Roub{\'{\i}}{\v{c}}ek.
\newblock Geophysical models of heat and fluid flow in damageable poro-elastic
  continua.
\newblock {\em Cont. Mech. Thermodyn.}, 29:625--646, 2017.

\bibitem{RoPaMa13QACV}
T.~Roub\'{\i}\v{c}ek, C.~G. Panagiotopoulos, and V.~Manti\v{c}.
\newblock Quasistatic adhesive contact of visco-elastic bodies and its
  numerical treatment for very small viscosity.
\newblock {\em Zeitschrift angew. Math. Mech.}, 93:823--840, 2013.

\bibitem{RoSoVo13MRLF}
T.~Roub\'{\i}\v{c}ek, O.~Sou\v{c}ek, and R.~Vodi\v{c}ka.
\newblock A model of rupturing lithospheric faults with re-occurring
  earthquakes.
\newblock {\em SIAM J. Appl. Math.}, 73:1460--1488, 2013.

\bibitem{Secc13VOC}
R.~A. Secco.
\newblock Viscosity of the outer core.
\newblock In T.~Ahrens, editor, {\em Mineral Physics \& Crystallography: A
  Handbook of Physical Constants}, pages 218--226. Willey, 2013.

\bibitem{SmyPal07VEOC}
D.~E. Smylie and A.~Palmer.
\newblock Viscosity of {E}arth's outer core.
\newblock Preprint arXiv:0709.3333, 2007.

\bibitem{ToCaMa09SSLV}
N.~Tosi, O.~\v{C}adek, and Z.~Martinec.
\newblock Subducted slabs and lateral viscosity variations: effects on the
  long-wavelength geoid.
\newblock {\em Geophys. J. Int.}, 179:813--826, 2009.

\bibitem{TAKS13EEEE}
V.~Tsai, J.-P. Ampuero, H.~Kanamori, and D.~Stevenson.
\newblock Estimating the effect of {E}arth elasticity and variable water
  density on tsunami speeds.
\newblock {\em Geophysical Res. Letters}, 40:492--496, 2013.

\bibitem{DKVD98VLIP}
G.~A.~D. Wijs, G.~Kresse, L.~Vo\v{c}adlo, D.~Dobson, D.~Alf\'e, M.~J. Gillan,
  and G.~D. Price.
\newblock The viscosity of liquid iron at the physical conditions of the
  {E}arth's core.
\newblock {\em Nature}, 392 (6678):805--807, 1998.

\bibitem{WooDeu09TOEF}
J.~H. Woodhouse and A.~Deuss.
\newblock Theory and observations - {E}arth{'}s free oscillations.
\newblock In B.~Romanowicz and A.~Dziewonski, editors, {\em Seismology and
  Structure of the Earth: Treatise on Geophysics}, volume~1, chapter 1.02,
  pages 31--65. Elsevier, 2009.

\end{thebibliography}
\end{document}